\documentclass[11pt,reqno]{amsproc}
\usepackage{enumerate}
\usepackage{multirow}
\usepackage{subfigure}
\usepackage{hyperref}
\usepackage{graphicx}
\usepackage{fullpage}
\usepackage{rotating}
\usepackage{color}

\usepackage{sidecap}
\usepackage{fullpage}
\usepackage{rotating}
\usepackage{epstopdf}
\usepackage{psfrag}
\usepackage{wrapfig}
\usepackage[small,bf]{caption}

\newtheorem{theorem}{Theorem}[section]
\newtheorem{remark}[theorem]{Remark}

\title{A framework for coupling flow and deformation of the porous solid}
\author{D.~Z.~Turner}
\address{Dr.~Daniel Z.~Turner, Department of Civil Engineering, 
University of Stellenbosch, Stellenbosch, South Africa. 
TEL: +27-21-808-4434} \email{dzturner@sun.ac.za}
\author{K.~B.~Nakshatrala}
\address{Dr.~Kalyana Babu Nakshatrala, Department of Civil 
and Environmental Engineering, University of Houston, Houston, 
Texas - 77204. TEL:+1-713-743-4418} \email{knakshatrala@uh.edu}
\author{M.~J.~Martinez}
\address{Dr.~Mario J.~Martinez, Thermal and Fluid Processes 
Department, Sandia National Laboratories, Albuquerque, New 
Mexico - 87185-0836. TEL:+1-505-844-8729} 
\email{mjmarti@sandia.gov}

\begin{document} 

\begin{abstract}
In this paper, we consider the flow of an incompressible fluid in a 
\emph{deformable} porous solid. We present a mathematical model 
using the framework offered by the theory of interacting continua. 
In its most general form, this framework provides a mechanism for capturing
multiphase flow, deformation, chemical reactions and thermal processes,
as well as interactions between the various physics in a conveniently implemented
fashion. To simplify the presentation of the framework, results are presented for
a particular model than can be seen as an extension of Darcy's equation (which 
assumes that the porous solid is rigid) that takes into account elastic deformation 
of the porous solid. The model also considers the effect of deformation 
on porosity. We show that using this model one can recover identical
results as in the framework proposed by Biot and Terzaghi. Some salient features of the framework are as 
follows: (a) It is a consistent mixture theory model, and adheres to 
the laws and principles of continuum thermodynamics, (b) the model 
is capable of simulating various important phenomena like consolidation 
and surface subsidence, and (c) the model is amenable to several 
extensions. 
We also present numerical coupling algorithms to obtain coupled 
flow-deformation response. Several representative numerical 
examples are presented to illustrate the capability of the mathematical 
model and the performance of the computational framework.
\end{abstract}


\maketitle

\section{MOTIVATION AND INTRODUCTION}
\textsc{Higher} oil prices and increasing awareness of the 
environmental impact of carbon pollution have motivated 
substantial interest in carbon-dioxide capture and storage 
(CCS). There is a growing consensus in both policy circles 
and in the energy industry that within the next few years, the 
US federal government will adopt some form of regulation for 
$\mathrm{CO}_2$ emissions \cite{Leach_Mason_Veld_NBER_2009}. 
At the same time, it is widely believed that much of the energy supply 
over the coming decades will continue to come from fossil 
fuels. Many analysts believe that the only way to reconcile 
the anticipated growth in the use of fossil fuels with anticipated 
limits on $\mathrm{CO}_2$ emissions is through the development 
and deployment of carbon-dioxide capture and sequestration. 

Of the proposed techniques for abating anthropogenic 
carbon-dioxide, storage in deep geologic formations 
is the only method that has gained widespread acceptance 
for its feasibility \cite{Macfarlance_Elements_2007_v3_p165,
  Schrag_Elements_2007_v3_p171}. If geological carbon-dioxide 
sequestration is implemented on the scale needed to make noticeable 
reductions in atmospheric $\mathrm{CO}_2$, a billion metric tons or 
more must be sequestered annually -- a 250-fold increase over the 
amount sequestered today. Large sedimentary basins are considered 
best suited to sequester such large volumes of $\mathrm{CO}_2$ as 
they have tremendous pore volume and connectivity and they are 
widely distributed \cite{Bachu_EG_2003_v44_p277,
  Benson_Cole_Elements_2008_v4_p325}. 
  
The development of energy systems like geological carbon-dioxide 
sequestration requires the understanding and predictive simulations 
of complex processes in earth systems. In particular, securing a 
large volume of carbon-dioxide will require a solid scientific foundation 
defining the coupled hydrologic-geochemical-geomechanical processes 
that govern the long term fate of  $\mathrm{CO}_2$. 
This path is becoming more dependent on modeling coupled geomechanical, thermal, fluid and chemical processes and the response of natural and engineered systems. 
Similarly, these models can be utilized to simulate oil and gas reservoirs, 
advanced recovery from tar sands and shales, underground storage 
of hydrogen, natural gas and oil, and evaluating aquifers for new water 
resources. Significant progress in these areas requires 
new coupled modeling approaches. 

\begin{SCfigure}
  \includegraphics[width=0.37\textwidth]{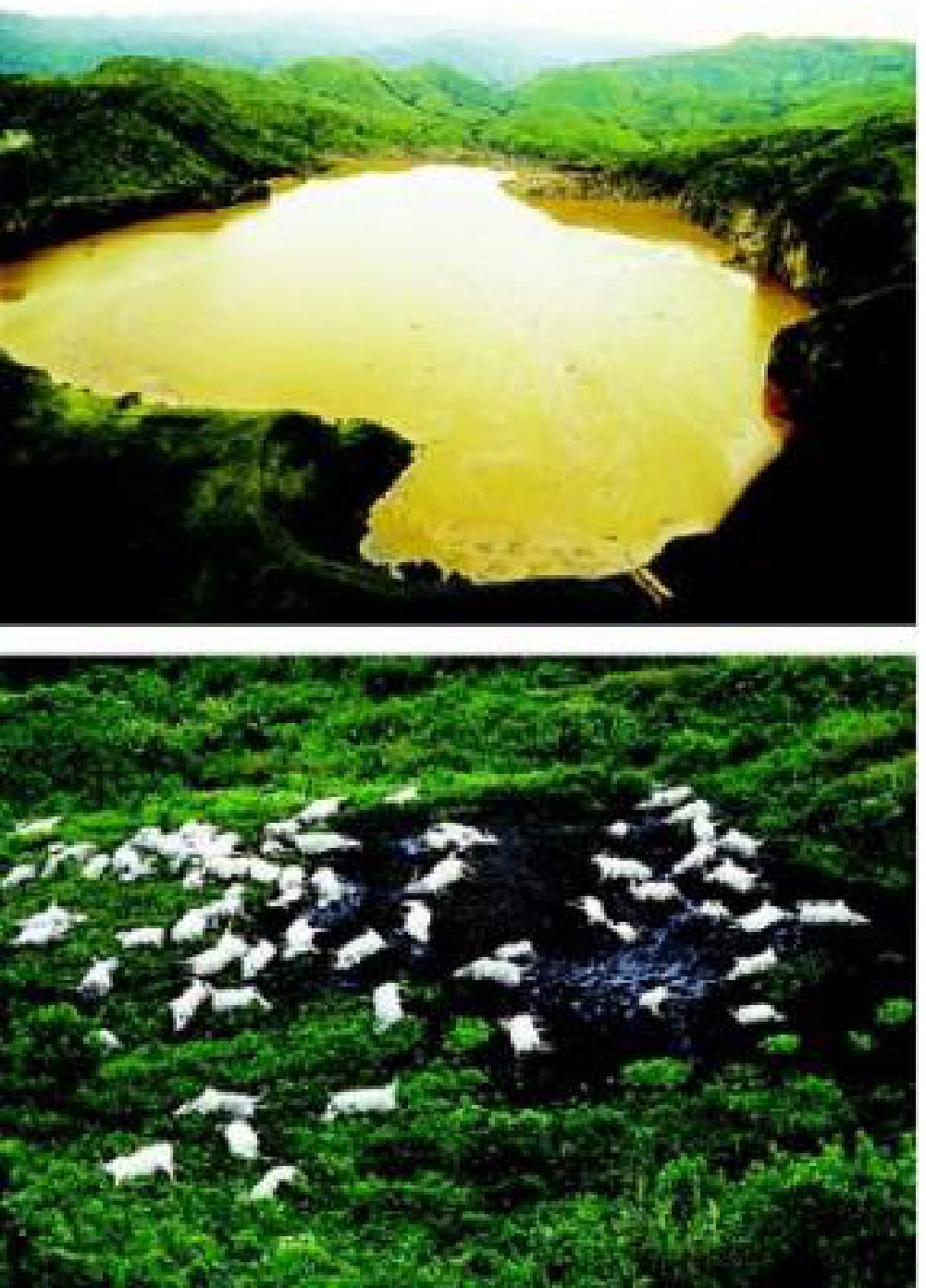}
  \caption{Lake Nyos is a deep crater lake in Cameroon (top figure). 
  In 1986, huge amounts of carbon-dioxide gas were released from 
  the lake possibly due to seismic tremors. A mix of carbon-dioxide 
  and water erupted 120 meters above the lake with a speed of 
  roughly 100 kilometers per hour. This displaced oxygen resulting 
  in asphyxiation of more than 1,700 people and 3,000 cattle as 
  far as 23 kilometers from the lake (bottom figure). 
  [Source: http://www.waterencyclopedia.com/] 
  \label{Fig:ECRP_Lake_Nyos}}
\end{SCfigure}

Unfortunately, there are still a number of unanswered scientific 
questions regarding deep geological sequestration that critically 
affect the efficacy and safety of this method. Although a number 
of studies have been conducted on geochemical interactions within 
the storage reservoir, and on the flow characteristics of the 
carbon-dioxide plume, relatively little effort has been devoted 
towards the reservoir's structural integrity. The seal created 
by the cap rock is one of the primary mechanisms that prevents 
injected carbon-dioxide from escaping back into the atmosphere. 
Failure of this seal can release large quantities of carbon-dioxide, 
which will have dire consequences. 
For example, in 1986 more than 1700 people were killed by the sudden 
release of geologic carbon-dioxide from lake Nyos in Cameroon (see 
Figure~\ref{Fig:ECRP_Lake_Nyos}). Although, in the case of lake Nyos, 
the carbon-dioxide was sequestered naturally, one can learn valuable 
lessons from this disaster. What happened in lake Nyos is similar to 
opening a shaken bottle of carbonated soda. As long as the cap is on, 
the carbon-dioxide gas stays dissolved under pressure. But when the 
cap is removed, the bubbles (and the soda) rapidly flow out of the 
bottle. The lake Nyos incident clearly highlights the importance of 
the structural integrity of the cap rock, and an immediate need for 
a systematic study along the lines presented in this paper. 

It is important to note that very few parallel processing commercial and/or 
research software tools exist for simulating complex processes such as 
coupled multiphase flow with chemical transport and geomechanics. 
Current computational limitations place significant restrictions on realistic 
problems that can be solved. Predictive computational simulation of the 
coupled-physics associated with geosystems is a critical enabling technology 
for their optimal management. A major stumbling block to high-fidelity 
modeling and simulation of geosystems is the lack of viable, robust, and 
efficient computational technologies for describing multiphase, multicomponent 
chemically reactive fluid mixtures in heterogeneous deformable geologic media. 
The present paper aims at advancing mathematical and numerical modeling, 
and predictive simulation of flow in deformable porous media under high 
pressures. 

\subsection{Main contributions of this paper}
Some of the main contributions of this paper are as follows:
\begin{enumerate}[(a)]
\item We have presented a mathematical model based on the theory 
of interacting continua which is capable of capturing surface 
subsidence and consolidation of soils. The model is fully coupled, 
and the interaction between the fluid and the porous solid is 
modeled through drag-like term and deformation-dependent porosity.  
\item We have presented various numerical coupling algorithms, 
  which can be used to obtain the coupled response.  
\item We have presented various representative numerical examples 
  illustrating the predictive capabilities of the proposed mathematical 
  model, and the numerical performance of the coupling algorithms. 
\end{enumerate}

\subsection{Organization of the paper}
The remainder of this paper is organized as follows. In Section 
\ref{Sec:Coupled_TIC} we give a brief review of the theory of 
interacting continua. In Section \ref{Sec:Coupled_Model} we 
present a model for the flow of an incompressible fluid in deformable 
porous solid. In Section \ref{Sec:Coupled_Algorithms} we shall 
present coupling algorithms that will be employed to obtained 
the coupled flow-deformation response. Several representative 
numerical results will be presented in Section \ref{Sec:Coupled_NR}, 
and conclusions are drawn in Section \ref{Sec:Coupled_Conclusions}.

\section{ELEMENTS OF THE THEORY OF INTERACTING CONTINUA}
\label{Sec:Coupled_TIC}
We shall employ the mathematical framework offered by the theory of 
interacting continua (TIC), which is sometimes referred to as the 
theory of mixtures. The basic assumption of the theory of interacting 
continua is that the constituents can be homogenized and assumed to 
co-occupy the domain occupied by the mixture. The theory of interacting 
continua traces its origins to the works of Darcy \cite{Darcy_1856} (also 
see its English translation by Patricia Bobeck \cite{Bobeck_Darcy}) and 
Fick \cite{Fick_1855_v94_p59}. Truesdell later gave the theory a firm 
mathematical footing (see \cite{Truesdell_Rendiconti_1957_v22_p33,
Truesdell_Rendiconti_1957_v22_p158,Truesdell_Toupin} and several 
appendices in Reference \cite{Truesdell_rational_thermodynamics}). 
We now provide a brief review of the basic equations of the theory of 
interacting continua. A more detailed treatment can be found in Atkin 
and Craine \cite{Atkin_Craine_QJMAM_1976_v29_p209}, Bowen \cite{Bowen}, 
and Bedford and Drumheller \cite{Bedford_Drumheller_IJES_1983_v21_p863}. 
One can also consult texts by Bowen \cite{Bowen_Porous_elasticity}, Coussy 
\cite{Coussy_poromechanics,Coussy_2010}, de Boer \cite{Boer_TPM}, Rajagopal 
and Tao \cite{Rajagopal_Tao}, and Voyiadjis and Song \cite{Voyiadjis_Song}. 

\begin{figure}[t]
  \psfrag{O}{$\delta V$}
  \psfrag{O1}{$\delta V^{(1)}$}
  \psfrag{O2}{$\delta V^{(2)}$}
  \psfrag{ON}{$\delta V^{(N)}$}
  \includegraphics[scale=0.55]{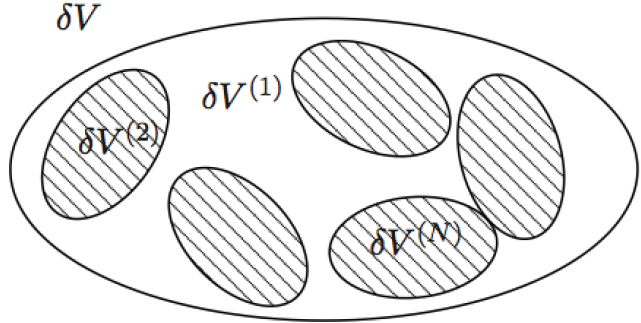}
  \caption{Elementary volume of the porous material $\delta V$, 
    and the volume occupied by the $i$-th constituent $\delta 
    V^{(i)}\; (i = 1, \cdots, N)$.}
\end{figure}

In the remainder of the paper, it should be noted that the 
usual summation convention on repeated indices will not be 
adopted. Consider a mixture of $N$ components. One of the 
components can be the porous solid, which will be the case 
when dealing with a deformable porous solid. A typical particle 
belonging to each constituent in the reference state is denoted 
$\boldsymbol{X}^{(i)} \; (i= 1, \cdots, N)$. At time $t$, these 
particles occupy the position $\boldsymbol{x}$. The motion and 
velocity of each constituent $(i = 1, \cdots, N)$ are defined 
through 
\begin{align}
  \boldsymbol{x} &= \chi^{(i)} (\boldsymbol{X}^{(i)},t) \\
  \boldsymbol{v}^{(i)} &= \frac{\partial \chi^{(i)}
    (\boldsymbol{X}^{(i)},t)}{\partial t}
\end{align}
We denote the gradient and divergence operators with respect to 
$\boldsymbol{x}$ through $\mathrm{grad}[\cdot]$ and $\mathrm{div}[\cdot]$, 
respectively. We now define the volume fractions and porosity by 
considering a representative volume element (RVE) with 
volume $\delta V$. The volume occupied by the $i$-th 
component in the RVE is denoted $\delta V^{(i)}$. The volume fraction of the 
$i$-th component is then defined as follows:
\begin{align}
  \eta^{(i)} := \frac{\delta V^{(i)}}{\delta V}
\end{align}
By definition, the volume fractions satisfy the following relationship:
\begin{align}
  \sum_{i=1}^{N} \eta^{(i)}(\boldsymbol{x}) = 1
\end{align}
The porosity of the porous medium is defined as follows:
\begin{align}
  \phi(\boldsymbol{x}) := 1 - \eta^{(s)}(\boldsymbol{x})
\end{align}
where $\eta^{(s)}(\boldsymbol{x})$ is the volume fraction of the 
porous solid. Let the mass of the $i$-th component in this RVE 
be $\delta m^{(i)}$. The true and bulk mass densities of the 
$i$-th component are, respectively, denoted by $\gamma^{(i)}$ 
and $\rho^{(i)}$. That is, 
\begin{align}
  \gamma^{(i)} &:= \frac{\delta m^{(i)}}{\delta V^{(i)}} \\
  \rho^{(i)} &:= \frac{\delta m^{(i)}}{\delta V} 
\end{align}
The bulk mass density of the mixture is defined as follows:
\begin{align}
  \rho(\boldsymbol{x},t) &:= \sum_{i = 1}^{N} \rho^{(i)}
  (\boldsymbol{x},t) 
\end{align}
The mixture velocity is defined as follows:
\begin{align}
   \boldsymbol{v}(\boldsymbol{x},t) &:= \frac{1}{\rho(\boldsymbol{x},t)} 
   \sum_{i=1}^{N} \rho^{(i)}(\boldsymbol{x},t) \boldsymbol{v}^{(i)}
   (\boldsymbol{x},t)
\end{align}
The diffusion velocity $\tilde{\boldsymbol{v}}^{(i)}(\boldsymbol{x},t)$ 
of the $i$-th constituent is defined as follows:
\begin{align}
  \tilde{\boldsymbol{v}}^{(i)}(\boldsymbol{x},t) = 
  \boldsymbol{v}^{(i)}(\boldsymbol{x},t) - \boldsymbol{v}
  (\boldsymbol{x},t)
\end{align}
Let $\alpha(\boldsymbol{x},t)$ be any quantity (scalar, vector, 
tensor) defined at the point $\boldsymbol{x}$ in the mixture at 
time $t$. We then define 
\begin{subequations}
  \begin{align}
    \frac{\mathrm{D}^{(i)} \alpha}{\mathrm{D} t} &:= \frac{\partial}
    {\partial t} \alpha(\boldsymbol{x},t) + \mathrm{grad} 
    \left[\alpha(\boldsymbol{x},t)\right]  \cdot \boldsymbol{v}^{(i)}
    (\boldsymbol{x},t) \\ 
    \frac{\mathrm{D} \alpha}{\mathrm{D} t} &:= \frac{\partial}{\partial t} 
    \alpha(\boldsymbol{x},t) + \mathrm{grad} \left[\alpha(\boldsymbol{x},t)
    \right] \cdot \boldsymbol{v}(\boldsymbol{x},t)
  \end{align}
\end{subequations}
The velocity gradient for the $i$-th constituent $\boldsymbol{L}^{(i)}$ 
and its symmetric part $\boldsymbol{D}^{(i)}$ are, respectively, defined 
through 
\begin{align}
  \boldsymbol{L}^{(i)} &:= \mathrm{grad}[\boldsymbol{v}^{(i)}] \\ 
  \boldsymbol{D}^{(i)} &:= \frac{1}{2} \left(\boldsymbol{L}^{(i)} + 
  {\boldsymbol{L}^{(i)}}^{T} \right)
\end{align}
We assume the existence of a partial traction vector 
$\boldsymbol{t}^{(i)}$ and a partial stress tensor $\boldsymbol{T}^{(i)}$ 
associated with each constituent of the mixture such that
\begin{align}
  \boldsymbol{t}^{(i)} = (\boldsymbol{T}^{(i)})^{T} \boldsymbol{n}_{S} 
\end{align}
where $\boldsymbol{n}_{S}$ is the normal to the surface $S$. We 
define the total traction and total stress through 
\begin{align}
  \boldsymbol{t} &= \sum_{i = 1}^{N} \boldsymbol{t}^{(i)} \\
  \boldsymbol{T} &= \sum_{i = 1}^{N} \boldsymbol{T}^{(i)}
\end{align}
so that we have 
\begin{align}
  \boldsymbol{t} = \boldsymbol{T}^{T} \boldsymbol{n}_{S}
\end{align} 

\subsection{Balance laws for constituents and the mixture}
The balance of mass for the $i$-th constituent takes the following form: 
\begin{align}
  \frac{\partial \rho^{(i)}}{\partial t} + \mathrm{div}[\rho^{(i)} 
  \boldsymbol{v}^{(i)}] = m^{(i)}
\end{align}
where $m^{(i)}$ is the mass supply for the $i$-th constituent. 
It is noteworthy that the theory of interacting continua can 
take into account the possible inter-conversion of mass due 
to chemical reactions between the constituents. The balance 
of mass for the mixture as a whole warrants that 
\begin{align}
  \sum_{i=1}^{N} m^{(i)} = 0
\end{align}
The balance of linear momentum of individual constituents can 
be written as follows: 
\begin{align}
\label{Eqn:Origin_balance_of_momentum}
  \rho^{(i)} \frac{\mathrm{D}^{(i)} \boldsymbol{v}^{(i)}}{\mathrm{D} t} 
  = \mathrm{div}[\boldsymbol{T}^{(i)}]^T  + \rho^{(i)} \boldsymbol{b}^{(i)} 
  + \boldsymbol{i}^{(i)}
\end{align}
where $\boldsymbol{b}^{(i)}$ is the external body force that acts on the 
$i$-th constituent, and $\boldsymbol{i}^{(i)}$ represents the interactive 
forces acting on the $i$-th constituent. That is, $\boldsymbol{i}^{(i)}$ 
arises due to the forces exerted by the other constituents on the $i$-th 
constituent by virtue of their being forced to co-occupy the domain of the 
mixture. Constitutive relations need to be specified for these interaction 
forces and they in general depend on not just the $i$-th constituent but 
also on the other constituents. Also, it follows from Newton's third law 
that 
\begin{align}
  \label{Eqn:TIC_Newtons_third_law}
  \sum_{i=1}^{N} \left(\boldsymbol{i}^{(i)} + m^{(i)} 
    \boldsymbol{v}^{(i)} \right) = \boldsymbol{0}
\end{align}
In the case of a mixture, even in the absence of angular momentum 
supply, the individual partial stresses need not be symmetric. However, 
the total stress will be symmetric in the absence of angular momentum 
supply, i.e.,
\begin{align}
  \label{Eqn:Balance_of_angular_momentum}
  \sum_{i=1}^{N} \boldsymbol{T}^{(i)} = \left(\sum_{i=1}^{N} 
    \boldsymbol{T}^{(i)}\right)^{T}
\end{align}
We shall assume that the individual partial stresses are symmetric, 
and hence equation \eqref{Eqn:Balance_of_angular_momentum} will be 
satisfied automatically. That is, 
\begin{align}
  \boldsymbol{T}^{(i)} = \left(\boldsymbol{T}^{(i)}\right)^{T} 
  \quad \forall i = 1, \cdots, N
\end{align}
One can similarly write governing equations for the balance of 
energy and the second law of thermodynamics. Herein, we shall 
restrict the studies to isothermal processes, and hence these 
equations are either trivially satisfied or not coupled with the 
mechanical aspects of the problem.

\section{A MODEL FOR FLUID FLOW IN DEFORMABLE POROUS SOLID}
\label{Sec:Coupled_Model}
The above framework offered by the theory of interacting continua is 
quite general and can model a wide variety of interesting phenomena. 
Herein, however, we shall consider a simplified but representative 
model for fluid and solid, which will be a special case of the above 
comprehensive framework. This simplified model will be used in our 
study on the stability of coupling algorithms. It should be noted 
that even this simplified model can capture many interesting 
features in the flow of fluids in deformable solids as discussed 
later in this paper.

We do not consider chemical reactions, and hence can set 
\begin{align}
  m^{(i)} = 0 \quad \forall i = 1, \cdots, N
\end{align} 
The constitutive equation for the fluid is 
taken to be that of a perfect fluid. That is, 
\begin{align}
  \label{Eqn:Coupled_partial_stress_of_fluid}
  \boldsymbol{T}^{(f)} = -p^{(f)} \boldsymbol{I}
\end{align}
where $\boldsymbol{T}^{(f)}$ is the stress in the fluid, 
$\boldsymbol{I}$ is the second-order identity tensor, 
and $p^{(f)}$ is the pressure in the fluid. We take into 
account interactions at the fluid-solid interface by 
including the following \emph{drag-like} terms:
\begin{subequations}
  \begin{align}
    \boldsymbol{i}^{(f)} &= \alpha^{(fs)} \left(\boldsymbol{v}^{(f)} - 
      \boldsymbol{v}^{(s)}\right) \\
    \boldsymbol{i}^{(s)} &= \alpha^{(sf)} \left(\boldsymbol{v}^{(s)} - 
      \boldsymbol{v}^{(f)}\right)
  \end{align}
\end{subequations}
By Newton's third law \eqref{Eqn:TIC_Newtons_third_law}, 
we have $\alpha^{(fs)} = \alpha^{(sf)}$. The drag term models 
the friction between the fluid and the porous solid, and it 
does not consider the friction within the fluid. The drag 
coefficient $\alpha^{(fs)}$ can be written as follows: 
\begin{align}
  \alpha^{(fs)} = \frac{\mu^{(f)}}{k}
\end{align} 
where $\mu^{(f)}$ is the coefficient of viscosity of the 
fluid, and $k$ is permeability. In this paper, we shall 
allow the coefficient of viscosity is allowed to depend 
on the pressure, and is modeled using Barus formula 
\cite{Barus_AJS_1893_v45_p87,
Nakshatrala_Rajagopal_IJNMF_2011_v67_p342}: 
\begin{align}
  \mu^{(f)}\left(p^{(f)}\right) = \mu^{(f)}_0 \exp
  \left[\beta^{(f)} p^{(f)}\right] 
\end{align}
where $\beta^{(f)}$ has units of $\mathrm{Pa}^{-1}$. 
The values of $\mu^{(f)}_0$ and $\beta^{(f)}$ of 
various organic liquids can found in references 
\cite{Bridgman,Barus_AJS_1893_v45_p87}. It is 
straightforward to show that the balance of 
linear momentum, given in equation 
\eqref{Eqn:Origin_balance_of_momentum}, becomes  
\begin{align}
  \label{Eqn:Coupled_fluid_LM}
  \rho^{(f)} \left(\frac{\partial \boldsymbol{v}^{(f)}}{\partial t} 
  + \mathrm{grad}[\boldsymbol{v}^{(f)}] \boldsymbol{v}^{(f)} 
  \right) + \alpha^{(fs)} \left(\boldsymbol{v}^{(f)} - 
  \boldsymbol{v}^{(s)} \right) + \mathrm{grad}[p^{(f)}] 
  = \rho^{(f)} \boldsymbol{b}^{(f)}
\end{align}
\begin{remark}
  A simple model that takes into account the friction between 
  the layers of the fluid has been proposed by Brinkman 
  \cite{Brinkman_ASR_1947_vA1_p27}. The model is referred 
  to as Brinkman equation or sometimes Darcy-Stokes equation. 
  Under this model, the partial stress in the fluid takes the 
  following form:
  \begin{align}
    \boldsymbol{T}^{(f)} = -p^{(f)} \boldsymbol{I} + 
    2 \mu^{(f)} \boldsymbol{D}^{(f)}
  \end{align}
  In such a case, the above expression for balance of linear 
  momentum \eqref{Eqn:Coupled_fluid_LM} should be altered as 
  follows:
  \begin{align}
    \rho^{(f)} \left(\frac{\partial \boldsymbol{v}^{(f)}}{\partial t} 
  + \mathrm{grad}[\boldsymbol{v}^{(f)}] \boldsymbol{v}^{(f)} 
  \right) + \alpha^{(fs)} \left(\boldsymbol{v}^{(f)} - 
  \boldsymbol{v}^{(s)} \right) + \mathrm{grad}[p^{(f)}]  
  - \mathrm{div}[2 \mu^{(f)} \boldsymbol{D}^{(f)}] = 
  \rho^{(f)} \boldsymbol{b}^{(f)}
  \end{align}
\end{remark}
The fluid is assumed to be incompressible, and hence the 
balance of mass for the fluid can be written as follows:
\begin{align}
  \label{Eqn:Coupled_mass}
  \frac{\partial \phi}{\partial t} + \mathrm{div}
  \left[\phi \boldsymbol{v}^{(f)}\right] = 0
  \end{align}
where $\phi$ is the porosity of the porous media. It should be 
noted that the formulation can be easily extended to variable 
density flows. We shall model the subsurface rock as a linearized 
elastic material. That is,
\begin{align}
  \mathop{\mathrm{max}}_{\boldsymbol{x},t} \; 
  \big\|\mathrm{grad}[\boldsymbol{u}^{(s)}]
  \big\|_{*} \ll 1
\end{align}
where $\|\cdot\|_{*}$ denotes the trace norm 
\cite{Lax_functional_analysis}, which is 
defined as follows:
\begin{align}
  \|\boldsymbol{A}\|_{*} = \mathrm{tr} \left[
    \sqrt{\boldsymbol{A}^{T}\boldsymbol{A}}\right]
\end{align}
where $\boldsymbol{A}$ is a second-order tensor, and 
$\mathrm{tr}[\cdot]$ denotes the trace. The (linearized) 
strain in the porous solid is given by 
\begin{align}
  \boldsymbol{\epsilon}^{(s)} := \frac{1}{2} 
  \left(\mathrm{grad}[\boldsymbol{u}^{(s)}] + 
  \mathrm{grad}[\boldsymbol{u}^{(s)}]^{T} \right)
\end{align}
Note that the velocity of the solid is given by 
\begin{align}
  \boldsymbol{v}^{(s)} = \frac{\partial \boldsymbol{u}^{(s)}}{\partial t}
\end{align}
The governing equation for the balance of linear 
momentum for the porous solid takes the following 
form: 
\begin{align}
  \label{Eqn:Coupled_LE_LM}
  \rho^{(s)} \frac{\partial^2 \boldsymbol{u}^{(s)}}{\partial t^2} 
  &- \alpha^{(fs)} \left(\boldsymbol{v}^{(f)} - \boldsymbol{v}^{(s)}
  \right)- \mathrm{div}[\boldsymbol{T}^{(s)}] = \rho^{(s)} 
  \boldsymbol{b}^{(s)} 
\end{align}
where $\boldsymbol{u}^{(s)}$ is the displacement of the solid. 
Note that the interaction effect has opposite signs in the 
two balance equations \eqref{Eqn:Coupled_fluid_LM} and 
\eqref{Eqn:Coupled_LE_LM} in virtue of Newton's third law. 
The constitutive equation for the porous solid can be written 
as follows: 
\begin{align}
\label{Eqn:Coupled_stress_strain}
  \boldsymbol{T}^{(s)} = \lambda^{(s)} \mathrm{tr}
  [\boldsymbol{\epsilon}^{(s)}] \boldsymbol{I} + 2 
  \mu^{(s)} \boldsymbol{\epsilon}^{(s)}
\end{align}
where $\boldsymbol{T}^{(s)}$ is the stress in the porous solid, 
and $\lambda^{(s)}$ and $\mu^{(s)}$ are the Lam\'e parameters of 
the porous solid. A few remarks about the above mathematical model 
are in order.
\begin{remark}
  The partial stresses for the fluid (given by equation 
  \eqref{Eqn:Coupled_partial_stress_of_fluid}) and for the 
  solid (given by equation \eqref{Eqn:Coupled_stress_strain}) 
  are symmetric, which is in accordance with the assumption 
  discussed in the previous section. 
\end{remark}
\begin{remark}
  Equations \eqref{Eqn:Coupled_fluid_LM}, \eqref{Eqn:Coupled_mass} 
  and \eqref{Eqn:Coupled_LE_LM} are coupled and are in terms of 
  three unknowns: $p^{(f)}$, $\boldsymbol{v}^{(f)}$, and $\boldsymbol{u}^{(s)}$. 
  In some of the specialized models discussed below, the number 
  of unknowns can be reduced by substituting one or more of the 
  equations in the above system into the other equations.
\end{remark}

\subsection{Quasi-static models}
One popular approach to incorporate the influence of the fluid 
motion on the deformation of the solid and \emph{vice versa} is 
by modeling the porosity of the porous media to be dependent 
on the deformation gradient in the solid. Typically, in this 
scenario, a number of approximations are made which are 
necessary to model the balance of momentum for the fluid 
using Darcy's law. If we ignore inertia in both the fluid 
and the solid, we can drop the first term in equations 
\eqref{Eqn:Coupled_fluid_LM} and \eqref{Eqn:Coupled_LE_LM}. 
Similarly, if we assume that the solid velocity is much 
smaller than the fluid velocity, the interaction term 
(second term in equations \eqref{Eqn:Coupled_fluid_LM} 
and  \eqref{Eqn:Coupled_LE_LM}) may be written as 
$\alpha^{(f)} \boldsymbol{v}^{(f)}$. 

Darcy's equation can be substituted into the interaction term 
for the fluid velocity. Rather than replace the interaction, 
$-\alpha^{(fs)} \boldsymbol{v}^{(f)}$, with the gradient of the 
fluid pressure, $\mathrm{grad}[p^{(f)}]$, and the fluid body 
force, $-\rho^{(f)} \boldsymbol{b}^{(f)}$, we shall incorporate 
the influence of the fluid pressure in the solid stress tensor 
as
\begin{align}
  \label{Eqn:SolidStressPorePressure}
  \boldsymbol{T}^{(s)} = \boldsymbol{T}_e^{(s)} + \boldsymbol{T}_p^{(s)}
\end{align}
where $\boldsymbol{T}_p^{(s)} = -p^{(f)}\boldsymbol{I}$. 
Partially incorporating the interaction term due to the 
fluid pressure in the solid stress as above is possible 
due to the following relationship, $\mathrm{grad}[p^{(f)}] 
= -\mathrm{div}[\boldsymbol{T}_p^{(s)}]$. The fluid body 
force, $\rho^{(f)}\boldsymbol{b}^{(f)}$, must of course be 
taken into account.

\subsection{Steady-state response}
The governing equations for steady-state response of flow of 
an incompressible fluid in a deformable porous solid can be 
written as follows:
\begin{subequations}
  \begin{align}
    &\mbox{Governing equations for fluid} \nonumber \\
    \label{Eqn:Coupled_LM_fluid}
    &\alpha^{(fs)}(p^{(f)}) \boldsymbol{v}^{(f)} + 
    \mathrm{grad}[p^{(f)}] = \rho^{(f)} \boldsymbol{b}^{(f)} \\
    \label{Eqn:Coupled_MB_fluid}
    &\mathrm{div}[\phi \boldsymbol{v}^{(f)}] = 0 \\
    &\mbox{Governing equations for solid} \nonumber \\
    &\mathrm{div}[\boldsymbol{T}^{(s)}] + \alpha^{(fs)}(p^{(f)}) 
    \boldsymbol{v}^{(f)} + \rho^{(s)} \boldsymbol{b}^{(s)} = \boldsymbol{0} \\
    &\boldsymbol{T}^{(s)} = \lambda^{(s)} \mathrm{tr}[\boldsymbol{\epsilon}^{(s)}] 
    \boldsymbol{I} + 2 \mu^{(s)} \boldsymbol{\epsilon}^{(s)} \\
    &\boldsymbol{\epsilon}^{(s)} := \frac{1}{2} 
    \left(\mathrm{grad}[\boldsymbol{u}^{(s)}] + \mathrm{grad}
               [\boldsymbol{u}^{(s)}]^T\right) \\
    &\mbox{Deformation-dependent porosity} \nonumber \\
    \label{Eqn:Coupled_dd_porosity}
    &\phi = \frac{\phi_0}{1 + (1 - \phi_0) 
      \mathrm{tr}[\boldsymbol{\epsilon}^{(s)}]}
  \end{align}
\end{subequations}
where $\phi_0$ is the porosity when the porous solid is unstrained (i.e., 
$\mathrm{tr}[\boldsymbol{\epsilon}^{(s)}]=0$). The unknown field variables 
are: velocity of the fluid $\boldsymbol{v}^{(f)}(\boldsymbol{x})$, pressure 
in the fluid $p^{(f)}(\boldsymbol{x})$, displacement of the porous solid 
$\boldsymbol{u}^{(s)}(\boldsymbol{x})$, and porosity $\phi(\boldsymbol{x})$. 
The above model (given by equations 
\eqref{Eqn:Coupled_LM_fluid}--\eqref{Eqn:Coupled_dd_porosity}) is nonlinear 
and fully coupled in the sense that neither the fluid model nor the solid 
model can be solved independently. 

\begin{remark}
  In the literature, one can find many other models 
  for deformation-dependent porosity. For example, 
  the following expression has been proposed: 
  \begin{align}
    \phi = 1 - \frac{1 - \phi_0}{\exp[\mathrm{tr}
      [\boldsymbol{\epsilon}^{(s)}]]}
  \end{align}
  Another model, which is used for large deformations of 
  the porous solid, takes the following form: 
  \begin{align}
    \phi = C_{\phi} \left(1 - \frac{(1 - \phi_0)}{\mathrm{det}
      [\boldsymbol{F}^{(s)}]}\right)
  \end{align}
  where $\phi_0$ is the initial porosity, $\boldsymbol{F}^{(s)}$ 
  is the solid deformation gradient
  \begin{align}
    \boldsymbol{F}^{(s)} = \boldsymbol{I} + 
    \mathrm{grad}[\boldsymbol{u}^{(s)}]
  \end{align}
  and $C_{\phi}$ is a coefficient that may very in time, 
  and can be used to smooth jumps in the porosity in time. 
  In the case where the mass balance equation (equation 
  \eqref{Eqn:Coupled_MB_fluid}) is subcycled, large changes 
  in porosity can occur during synchronization updates. Using 
  a reference pressure, $p_{\mathrm{ref}}$, that is taken as the 
  pressure at the last synchronization update, the following 
  expression can be used to interpolate the porosity between 
  updates
  \begin{align}
    C_{\phi} = 1 + C_r(p - p_{\mathrm{ref}})
  \end{align}
  where $C_r$ is the rock compressibility constant.
\end{remark}

\section{NUMERICAL COUPLING ALGORITHMS}
\label{Sec:Coupled_Algorithms}
The coupled equations presented above can be solved in a number of ways, 
which will depend on how the overall problem is decomposed. The first method 
is the \emph{fully coupled} method, which solves governing equations for fluid, 
solid and porosity update simultaneously as one monolithic system. The variables 
in this case are solid displacement, fluid velocity, fluid pressure, and porosity. 
Under this method, it is typical to use the same time step for all subsystems. 
The second method is the \emph{lockstep} method (which is also referred to 
as Gauss-Seidel), which solves each subsystem individually until convergence 
in serial fashion.  The solution is then transferred to the other subsystems which 
then iterate until convergence. Under the lockstep method, the same time step is also used 
for all subsystems. In the third method, referred to as the \emph{subcycle} 
method, subcycling occurs in either the mass balance or solid deformation equation. 
For the problems presented below, subcycling was only used for the mass balance 
equation since the physics of the flow evolves much faster than the solid deformation. 
The last method considered in this paper is the \emph{Jacobi} method. Under this 
method, both equations are decoupled and converged individually, but fields are 
transferred after each iteration. Figure \ref{fig:CouplingAlgs} pictorially describes 
these coupling algorithms.

\section{REPRESENTATIVE NUMERICAL RESULTS}
\label{Sec:Coupled_NR}
In this section, using the mathematical model described in Section 
\ref{Sec:Coupled_Model}, we shall study the performance and relative 
efficiency of the coupling algorithms that are outlined in the previous 
section. Several canonical problems are solved. 

\subsection{Verification by a manufactured solution}
In this subsection, we shall use the method of manufactured solutions 
to verify the implementation of the coupling algorithms. The method 
of manufactured solutions assumes an exact solution and then derives 
the corresponding source terms that satisfy the governing equations. 
The (manufactured) exact solution for the coupled problem 
\eqref{Eqn:Coupled_LM_fluid}--\eqref{Eqn:Coupled_dd_porosity} 
takes the following form:
\begin{subequations}
  \begin{align}
    &u^{(s)}(x) = u_0^{(s)} \sin(\pi x), \; \epsilon^{(s)}(x) = 
    u_0^{(s)} \pi \cos(\pi x), \; T^{(s)} (x) = \left(\lambda^{(s)} 
    + 2 \mu ^{(s)} \right) u_0^{(s)} \pi \cos(\pi x) \\
    &v^{(f)}(x) = v_0^{(f)} \left(1 + (1 - \phi_0) u_0^{(s)} \pi 
    \cos(\pi x)\right), \; p^{(f)}(x) = -\alpha^{(fs)} v_0^{(f)} 
    \left(x + (1 - \phi_0) u_0^{(s)} \sin(\pi x)\right) + p_0^{(f)} \\
    &\phi(x) = \frac{\phi_0}{1 + (1 - \phi_0) \mathrm{tr}[\epsilon^{(s)}]} 
    = \frac{\phi_0}{1 + (1 - \phi_0) u_0^{(s)} \pi \cos(\pi x)}
  \end{align}
\end{subequations}  
The prescribed body forces for the fluid and the porous solid are as follows:
\begin{align}
  b^{(f)}(x) = 0, \quad b^{(s)}(x) = \frac{u_0^{(s)} \pi^2}{\rho^{(s)}} 
  \left(\lambda^{(s)} + 2 \mu^{(s)}\right) \sin(\pi x)
\end{align}
The constants used in this numerical experiment are listed in Table 
\ref{table:1DVerificationParams}.
\begin{table}[htb]
  \caption{Verification by a manufactured solution: Values used in 
    the numerical simulation. \label{table:1DVerificationParams}}
  \centering
  \begin{tabular}{c c c c}
    \hline
    Parameter & Value & Parameter & Value \\
    \hline
    $\rho^{(f)}$ & 1.0 & $\rho^{(s)}$ & 1.0 \\
    $v_0^{(f)} $  & 1.0 & $p_0^{(f)}$  & 1.0 \\
    $\lambda^{(s)} $  & 1.0 &$\mu^{(s)} $  & 0.5 \\
    $\alpha^{(fs)} $  & 1.0 & $\phi_0$  & 0.1 \\
    $u_0^{(s)}$  & 0.01 & \\ \hline
  \end{tabular}
\end{table}
To restrict the analysis of this problem to one-dimension, the 
boundary conditions are prescribed such that displacements in the 
$y$-direction are fixed and only displacements in the $x$-direction 
are considered. The problem domain consists of a unit square with 
a single element in the $y$-direction and 200 elements in the 
$x$-direction. On the left side of the domain the solid displacement 
is fixed. On the right side of the domain, the manufactured solution 
is applied for the solid displacements as a boundary condition.

\begin{sidewaysfigure}
\vspace*{15 cm}
	\centering
 \subfigure{\includegraphics[scale=0.5]{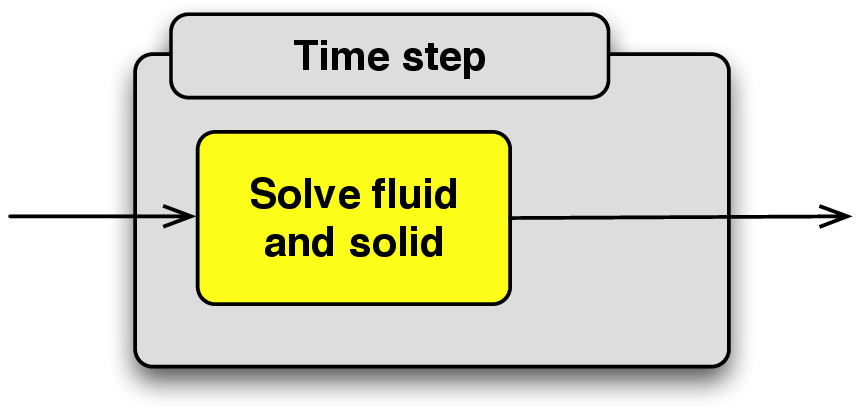}}
 \subfigure{\includegraphics[scale=0.5]{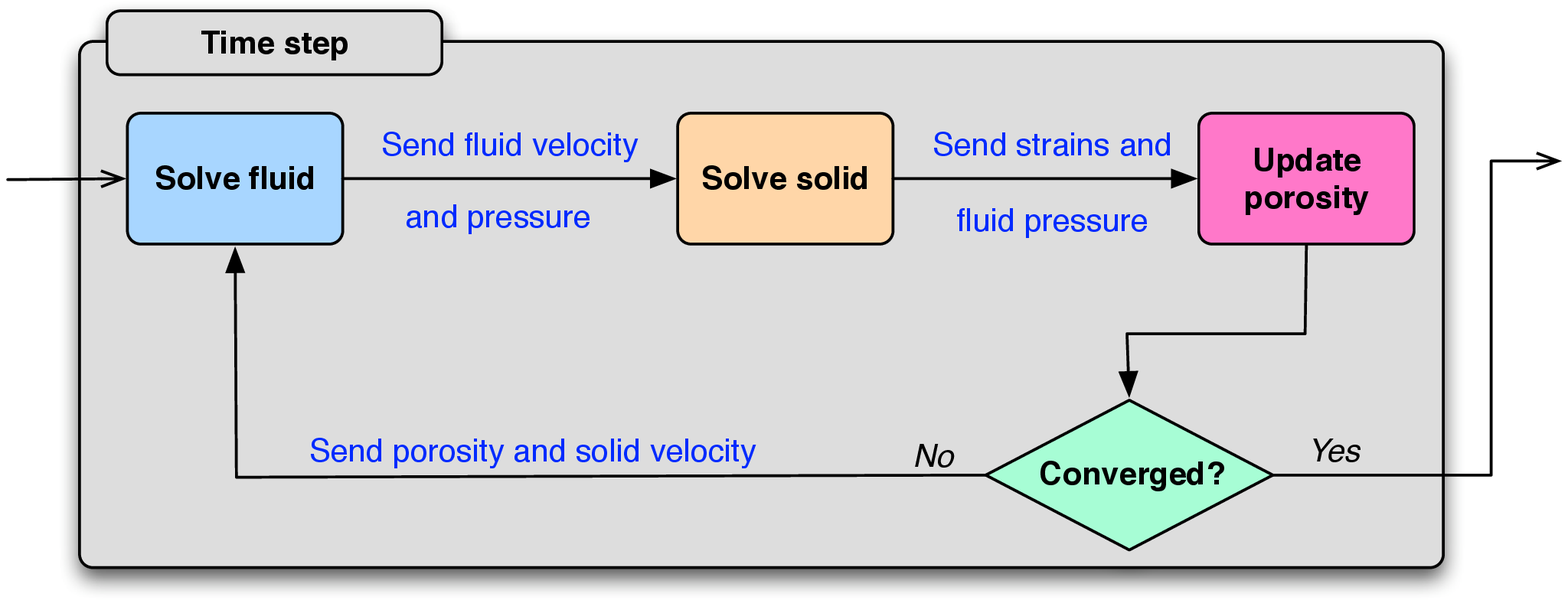}}
 \subfigure{\includegraphics[scale=0.5]{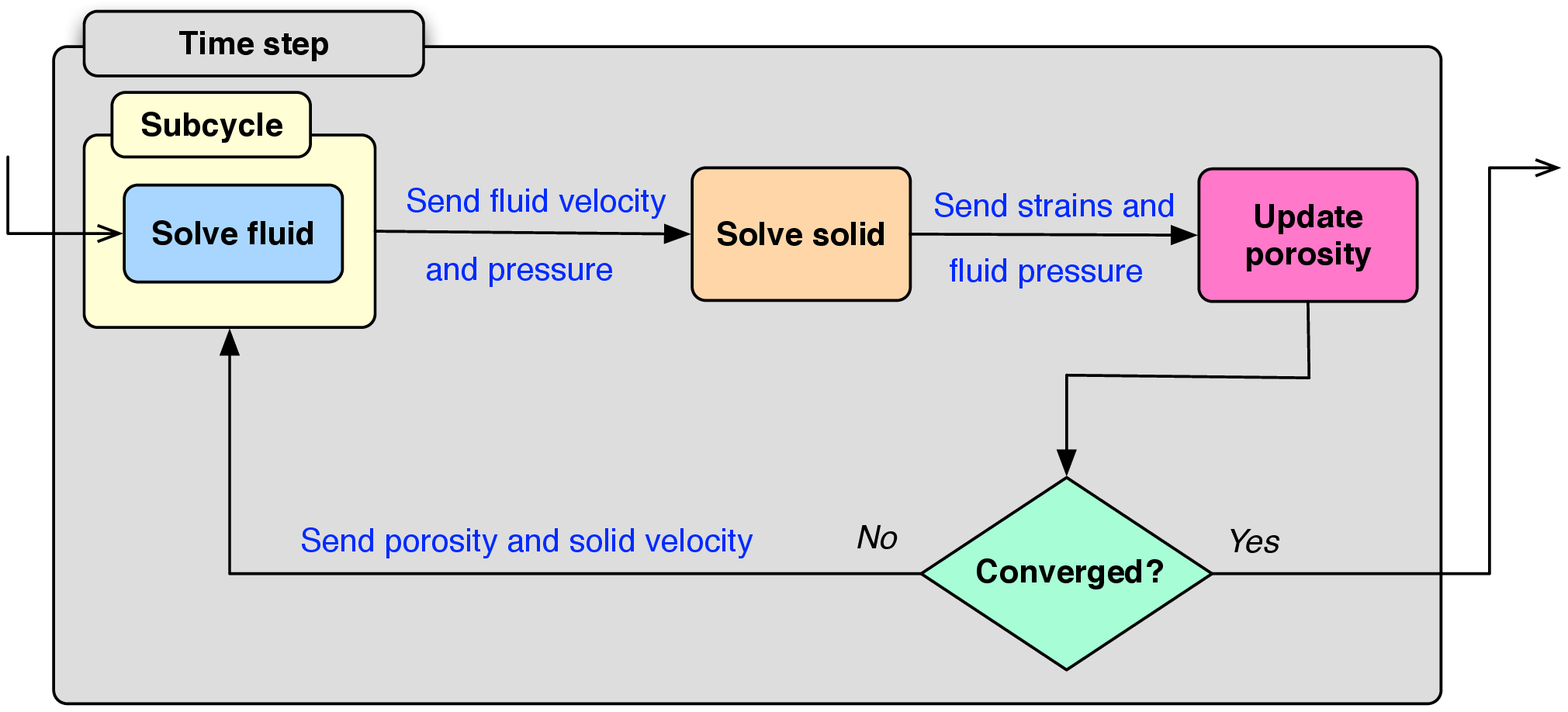}}
 \subfigure{\includegraphics[scale=0.5]{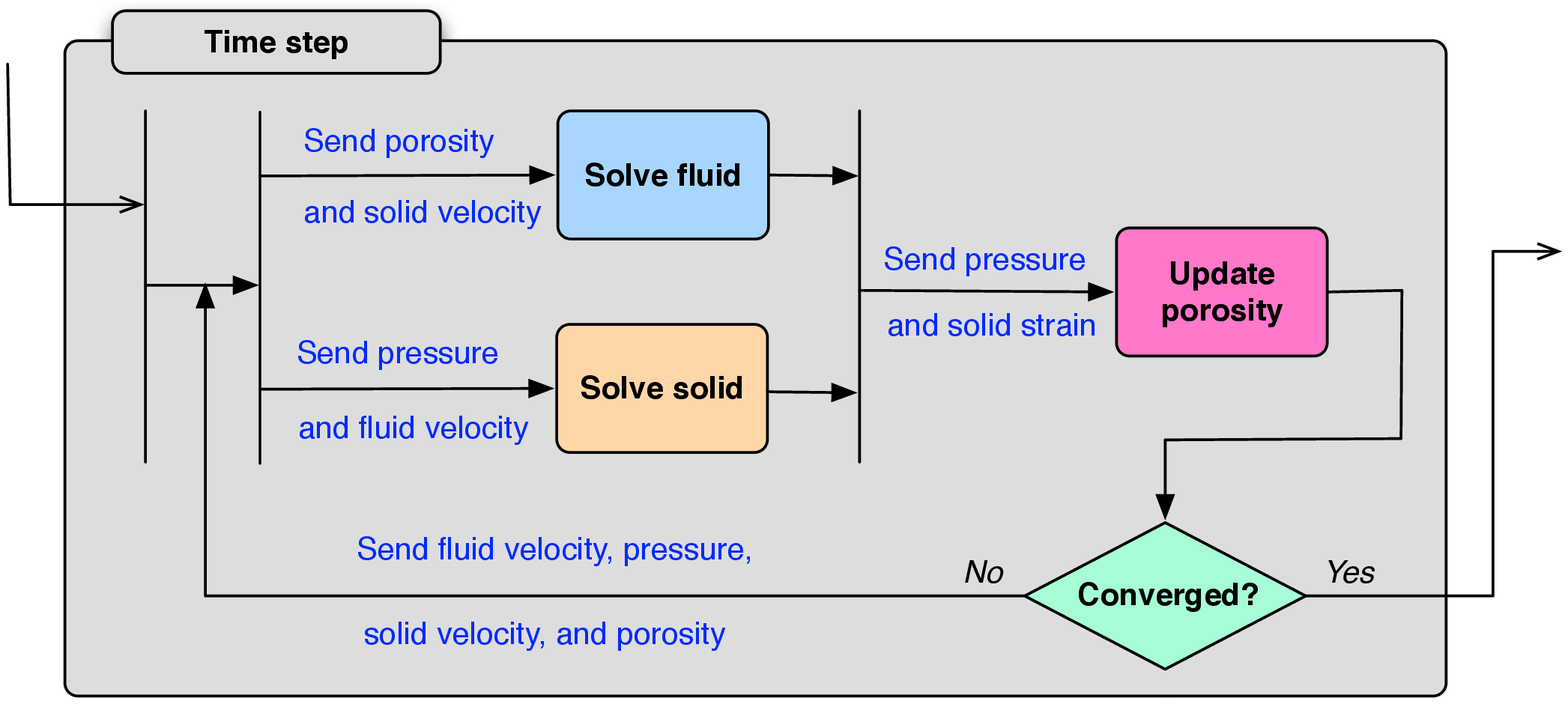}}
	\caption{Numerical coupling algorithms: fully coupled (top left), 
	lockstep (top right), subcycling (bottom left), and Jacobi (bottom right).} 
	\label{fig:CouplingAlgs}
\end{sidewaysfigure}

\begin{figure}
  \centering
  \subfigure{\includegraphics[scale=0.5]{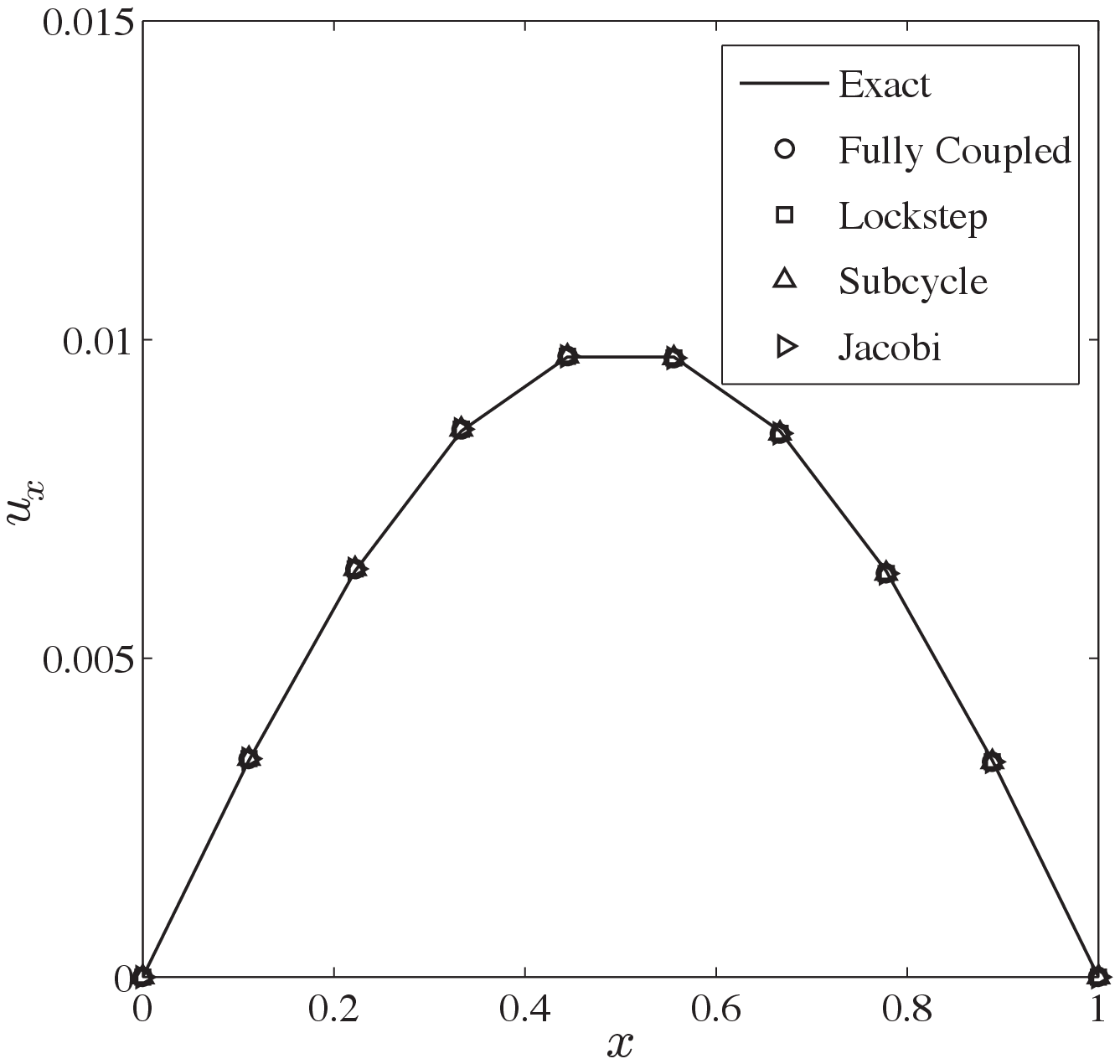}}
  \subfigure{\includegraphics[scale=0.5]{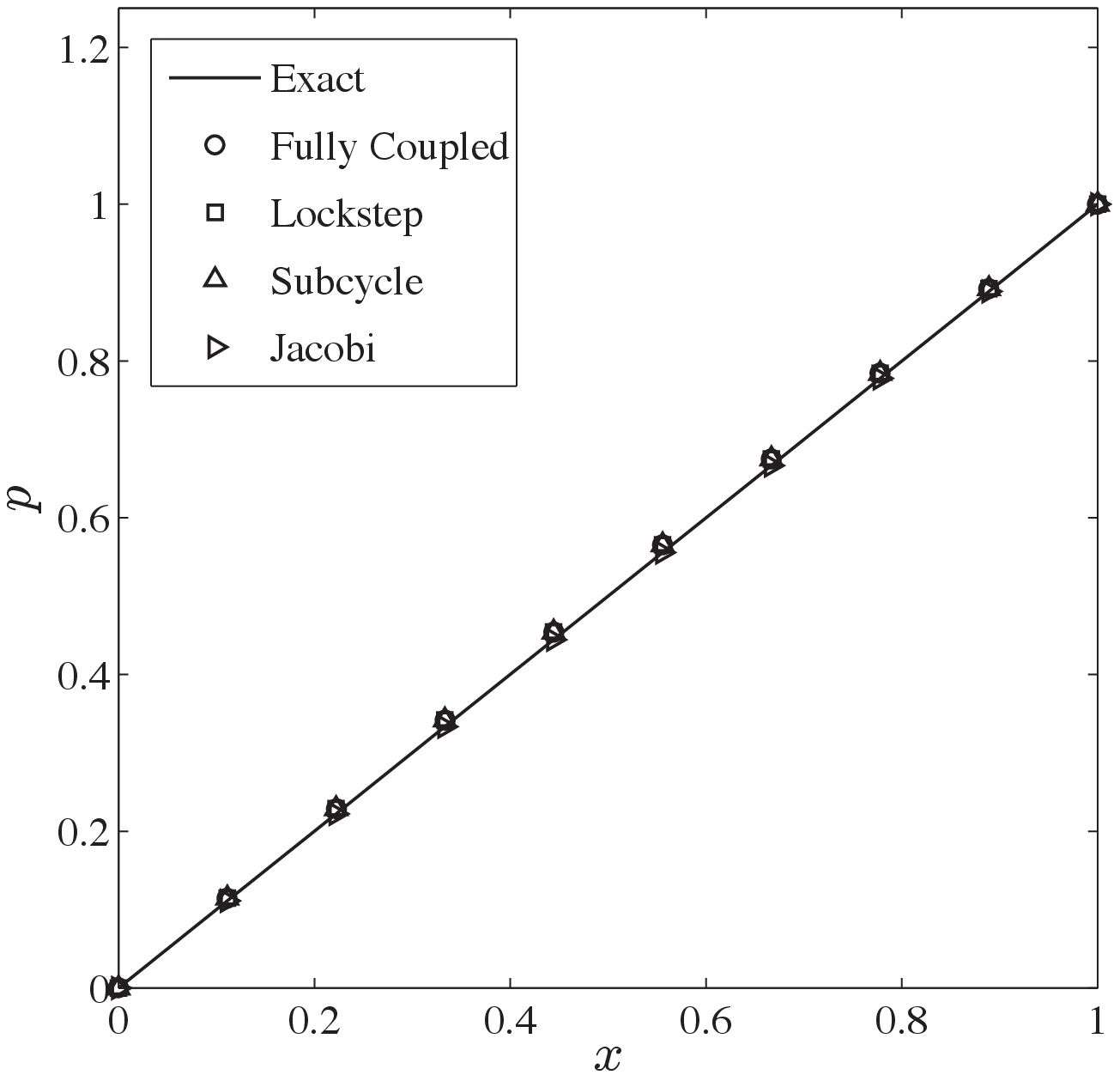}}
  \subfigure{\includegraphics[scale=0.5]{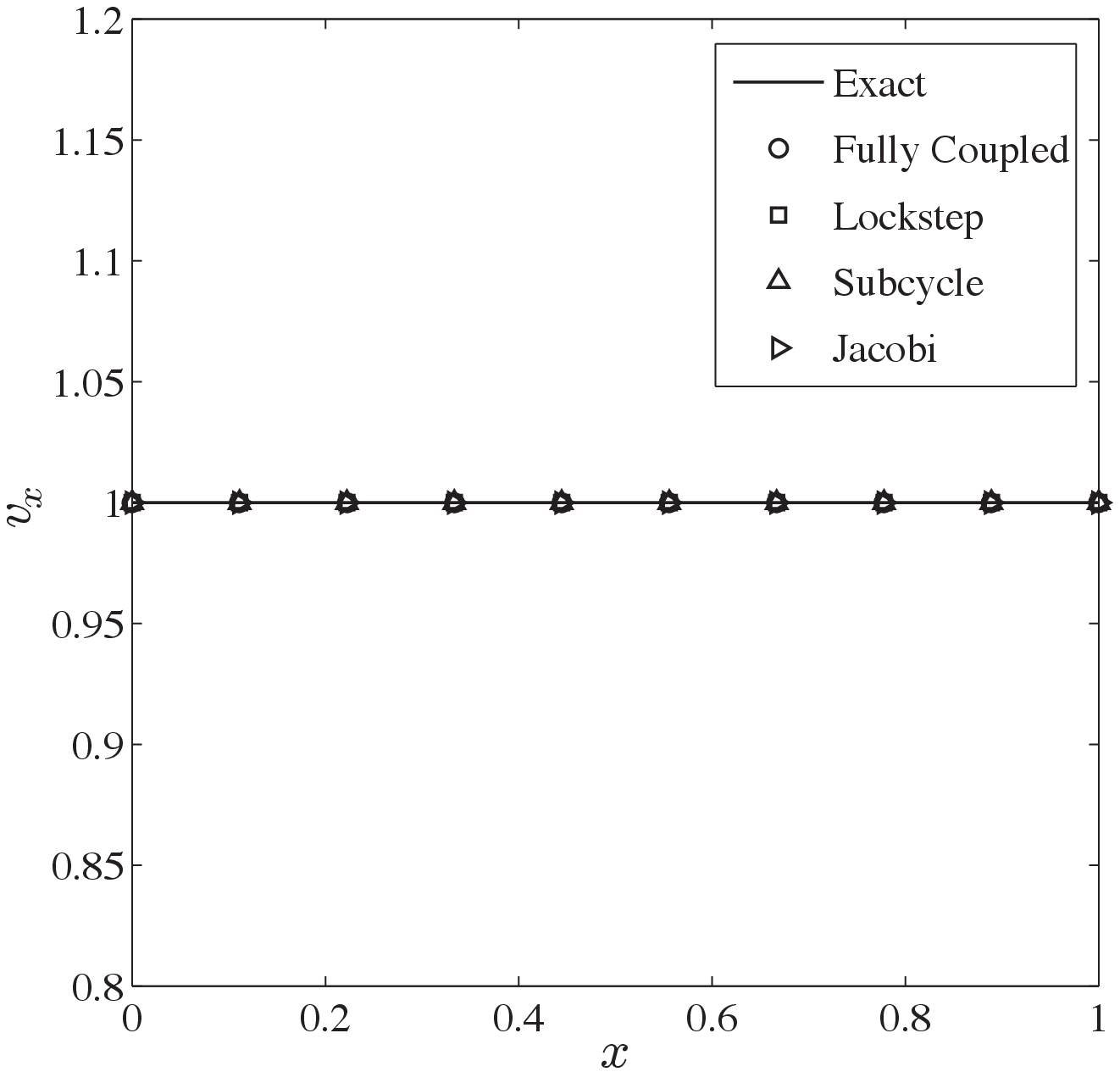}}
  \subfigure{\includegraphics[scale=0.5]{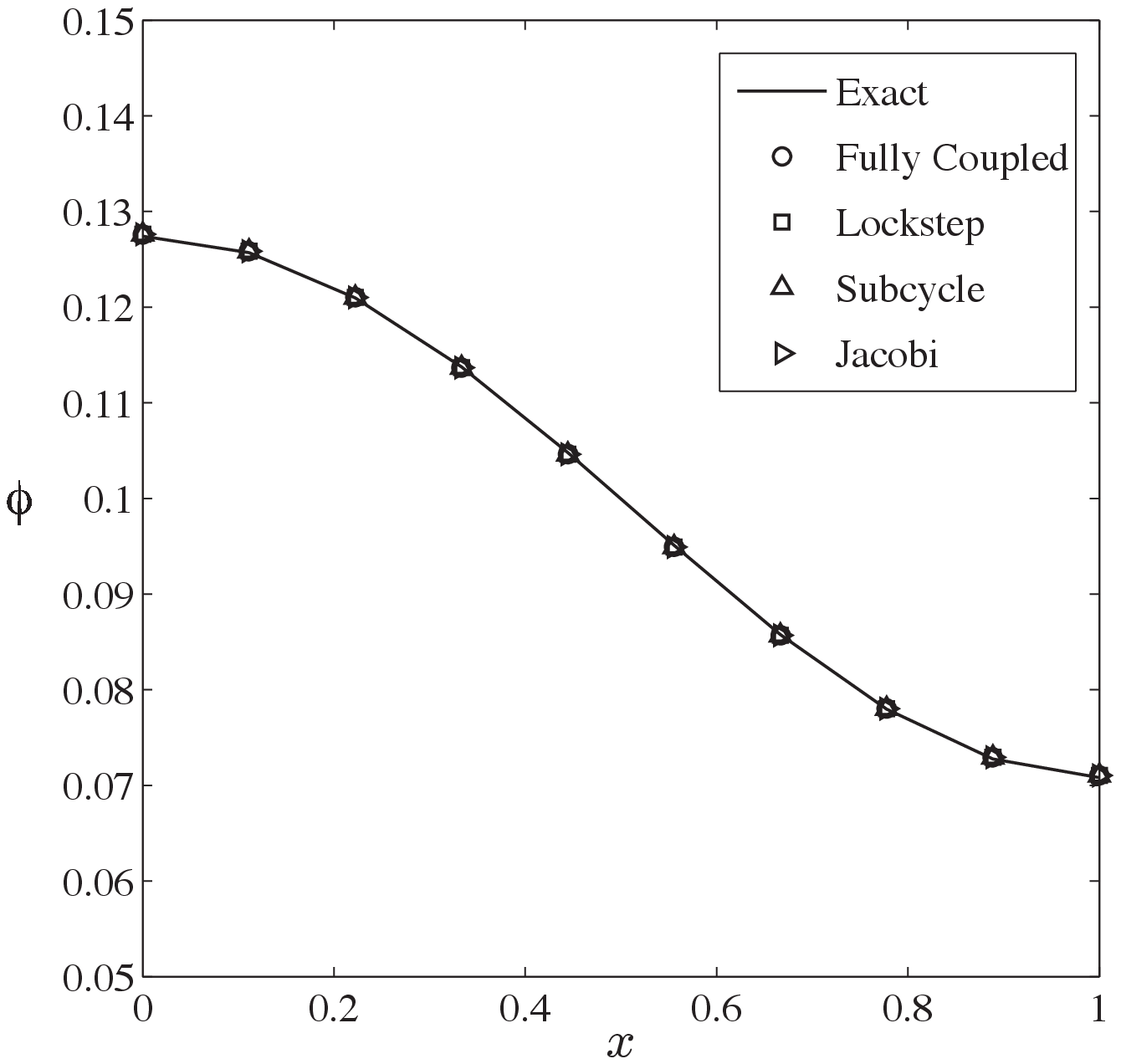}}
  \caption{Verification by a manufactured solution: numerical 
    solution vs. $x$ compared with the exact solution (top left) 
    solid displacement, $u_x$ (top right) pressure, $p$ (bottom left) 
    fluid velocity, $v_x$ (bottom right) porosity, $\phi$. The 
    \emph{lockstep} method was computed using a tolerance of 
    $\epsilon_{l} = 1 \times 10^{-9}$. The \emph{subcycle} method 
    was computed using one subcycle (equivalent to \emph{lockstep}) 
    and a tolerance of $\epsilon_{s} = 1 \times 10^{-9}$. The 
    \emph{Jacobi} method was computed using a tolerance of 
    $\epsilon_{j} = 1 \times 10^{-9}$.}\label{fig:ExactSolution}
\end{figure}

\begin{figure}
  \centering
  \includegraphics[scale=0.5]{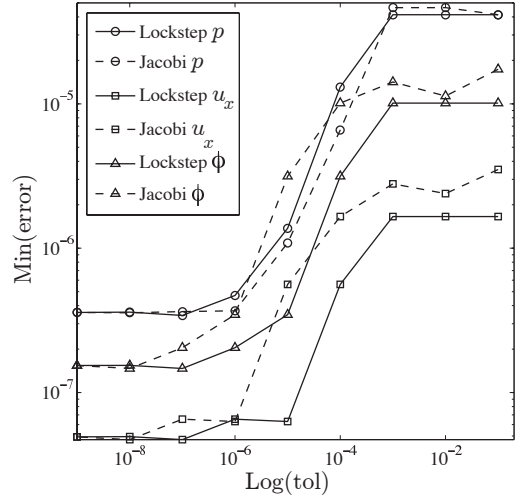}
  \caption{Verification by a manufactured solution: minimum 
    error vs. global convergence tolerance, $\epsilon_l$ or 
    $\epsilon_j$. The error is measured as the $L_2$-norm 
    of the computed solution minus the exact solution and 
    the minimum is taken over the total solution time.} 
  \label{fig:TolErrorLockstep1}
\end{figure}

\begin{figure}
  \centering
  \includegraphics[scale=0.5]{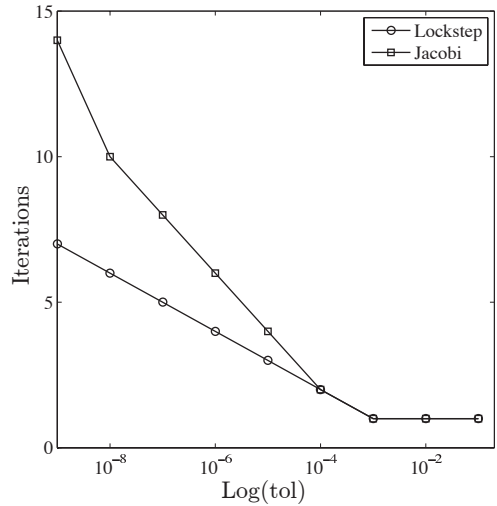}
  \caption{Verification by a manufactured solution: 
    number of iterations vs. global convergence 
    tolerance, $\epsilon_l$ or $\epsilon_j$. The 
    number of iterations is measured as the number 
    of fluid and solid solves per solution cycle
    before the global tolerance is met.}
  \label{fig:TolErrorLockstep2}
\end{figure}

The above outlined problem was solved by each of the coupling algorithms 
described in Figure \ref{fig:CouplingAlgs}. A comparison between the exact 
solution and the obtained numerical solution for the fluid pressure, solid 
displacements, and the porosity are shown in Figure \ref{fig:ExactSolution}. 
Figures \ref{fig:TolErrorLockstep1} and \ref{fig:TolErrorLockstep2} provide 
a measure of the efficiency of the loosely coupled algorithms (lockstep and 
Jacobi). The relative error of the solution is computed as the $L_2$ norm of 
the exact solution minus the computed solution. Notice that both the lockstep 
and Jacobi methods perform similarly in terms of the accuracy obtained for a 
given convergence tolerance, but the Jacobi method requires more iterations 
for a given level of accuracy.

\subsection{Terzaghi consolidation problem}
\begin{figure}[htb!]
  \centering
  \includegraphics[scale=1.0]{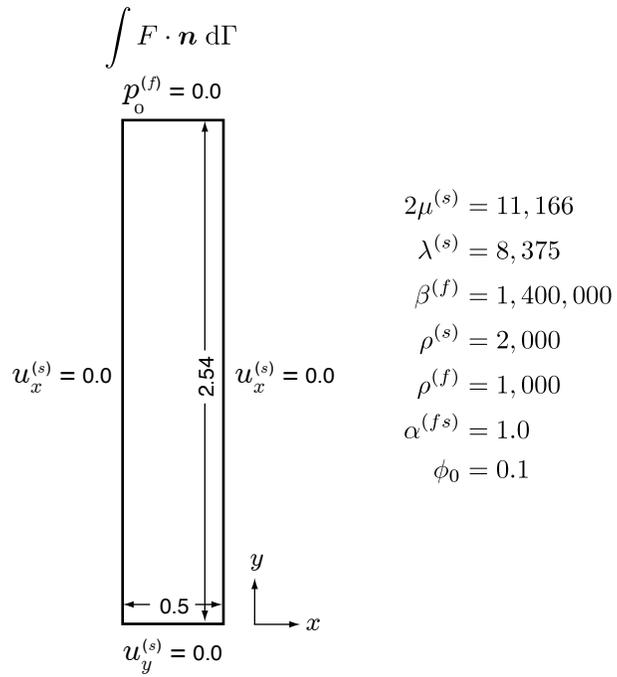}
  \caption{Terzaghi dynamic consolidation: problem geometry 
    and boundary conditions}
  \label{fig:TerzaghiDomain}
\end{figure}

To evaluate the performance of the \emph{quasi-static} model (equations \eqref{Eqn:Coupled_LM_fluid} to \eqref{Eqn:Coupled_dd_porosity}) for one-dimensional transient wave propagation in fluid saturated incompressible porous media we investigate the response of an infinite half-space to time dependent loading. This example problem, which has been studied by a number of researchers including \cite{deBoer, Diebels}, is represented computationally by the domain and parameters shown in Figure \ref{fig:TerzaghiDomain}, although the dimensions of the domain are insignificant regarding the resulting response. The top surface of the domain is subjected to a forcing function, $F$, given as:
\begin{align}
F = 100\left((1 - \mathrm{cos}(75t) \right)
\end{align}
The distributed force is applied in a weak fashion over the top surface. This surface is also treated as \emph{perfectly drained} meaning that the fluid is free to drain from the top surface and the pressure is equal to the ambient pressure. Along the sides and bottom of the domain the fluid is not allowed to drain. The displacement boundary conditions consist of fixing the $x$-displacement along the right and left boundary and the $y$-displacement on the bottom surface. This essentially allows for only one-dimensional consolidation in the $y$-direction.

\begin{figure}[htb!]
  \centering
  \includegraphics[scale=0.5]{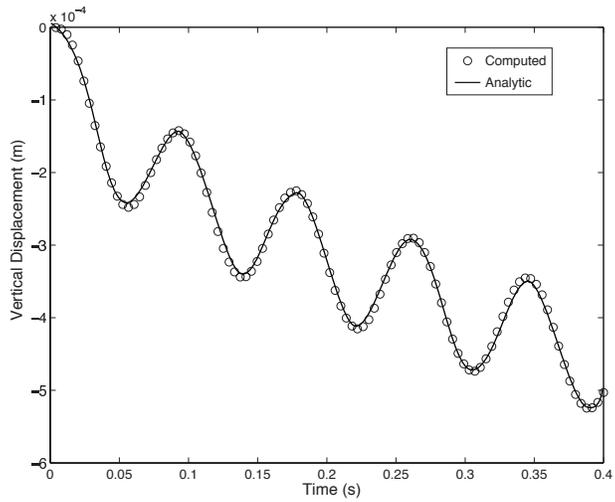}
  \caption{Terzaghi dynamic consolidation: computed solution 
    compared with the analytic solution given in \cite{deBoer}.}
  \label{fig:TerzaghiSolution}
\end{figure}

The analytic solution for the $y$-displacement of the top surface of the domain is given as:
\begin{align}
u^{(s)}_y(t) = \frac{-1.0}{\sqrt{a} (\lambda^{(s)} + 2\mu^{(s)})}\int_0^t F(t - \tau) e^{-\frac{b\tau}{2a}} I_0 \left( \frac{b\sqrt{\tau^2 - a z^2}}{2a} \right) U(\tau - \sqrt{a}z )  \; \mathrm{d}\tau
\end{align}
where $a = (n^{(s)2} \rho^{(f)} + n^{(f)2} \rho^{(s)})/c$, $b = n^{(f)2} k_c/c$, $c = (\lambda^{(s)} + 2\mu^{(s)}) n^{(f)2}$, $k_c = 1\times 10^6$, $I_0(z)$ is the modified Bessel function of zeroth order, and $U(t)$ is a unit Heaviside function. The computed solution using the proposed \emph{quasi-static} model is shown in Figure \ref{fig:TerzaghiSolution}. Notice that the model performs well even though the loading is nonlinear and the dynamic response of the solid is treated in a quasi-static fashion.

\subsection{Surface subsidence problem}
\begin{figure}[htb!]
  \centering
  \includegraphics[scale=0.9]{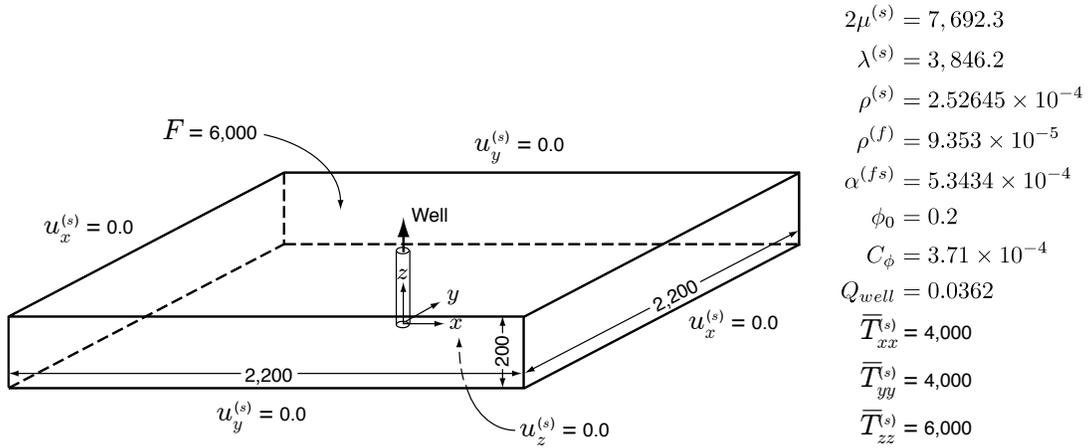}
  \caption{Surface subsidence: Problem geometry, parameters, and 
  boundary conditions.}
  \label{fig:DeanDomain}
\end{figure}

\begin{figure}[htb!]
  \centering
  \includegraphics[scale=0.5]{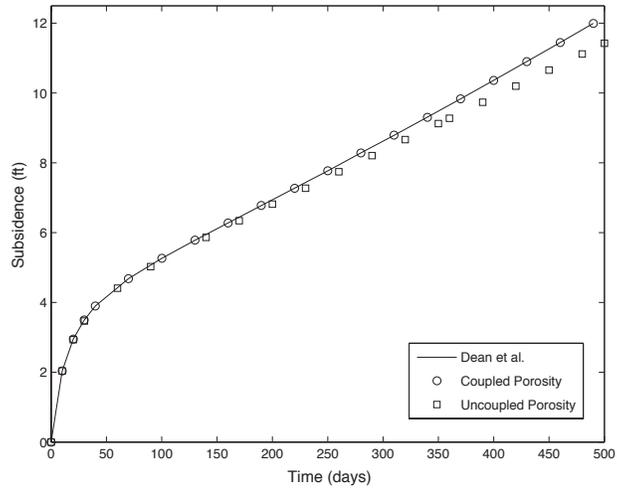}
  \caption{Surface subsidence: Comparison of the subsidence at the center of the top of the domain with published results in \cite{Dean} for both coupled and uncoupled porosity.}
  \label{fig:DeanCompare}
\end{figure}

\begin{figure}[htb!]
  \centering
  \includegraphics[scale=0.4]{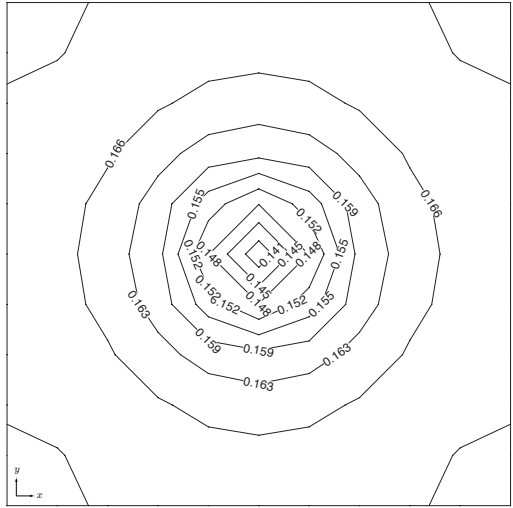}
  \caption{Surface subsidence: Contours of porosity for the top surface of the domain.}
  \label{fig:DeanPorosityContours}
\end{figure}

In this verification problem we consider the subsidence of a reservoir resulting from extraction of fluid via a centrally located well that has been studied extensively in \cite{Dean}. Figure \ref{fig:DeanDomain} shows the domain and boundary conditions used for this problem. In order to create confinement conditions similar to the subsurface environment an initial stress state was specified for the solid stress. The coupling between the solid stress and the fluid pore pressure is defined by equation \eqref{Eqn:SolidStressPorePressure} and the porosity model corresponds to equation \eqref{Eqn:Coupled_dd_porosity}. The \emph{quasi-static} model given in the previous examples was also used for the governing equations. Figure \ref{fig:DeanCompare} shows a comparison between the results published in \cite{Dean} and the computed results from the present study. Results are shown for both coupled and uncoupled porosity. Close correlation is achieved for the subsidence of the center of the top surface the domain for the coupled results. If the effect of the solid deformation is not taken into account in the porosity, the subsidence of the top surface is under-predicted by roughly 5\%. Figure \ref{fig:DeanPorosityContours} shows contours of porosity for the top surface of the domain. Notice that near the ejection well, the porosity is decreased by almost 30\% from the initial value.

\subsection{Two-phase immiscible water flood}
In the last example we present results for a water flood of 
the five spot well domain. The problem domain consists of 
two-dimensional square with injection wells at one of the 
sets of opposing corners and ejection wells at the other. 
Water is injected into the domain pushing oil out the 
ejection wells. The solution to this problem has been 
published for a number of numerical approaches including 
the results found in \cite{Helmig}. 

\subsubsection{Changes in porosity vs. change in permeability due to solid stress}
In the present study we consider the effects of both changes in \emph{porosity} due to stresses in the solid skeleton and resulting changes in \emph{permeability} due to damage. Conceptually, it is well established that solid deformation dependent porosity models behave like compressibility models. The otherwise incompressible fluid is able to be squeezed into pores that are expanding and vice versa. This is also apparent from the solid stress dependent contribution to the fluid equations showing up in the mass balance equation. In addition to the effect of pores enlarging and collapsing, one may wish to incorporate the increased permeability or interconnectedness of the void spaces due to solid damage. For example, in regions of high stress in the reservoir rock, cracks will form that allow fluid to flow with much less resistance. As a result, the flow magnitude and path of the fluid may change depending on the state of stress in the geologic formation. To explore differences between porosity coupled and permeability coupled models we investigate the two-phase water flood problem using both models separately. For the porosity coupled results we employ the model given in equation \eqref{Eqn:Coupled_dd_porosity}, for the permeability coupled results we consider a damage modified permeability model given as:
\begin{align}
\alpha^{(f)} = \alpha^{(f)}_0 \left(1.0 + \frac{\zeta || \boldsymbol{T}^{(s)} - \boldsymbol{T}_0^{(s)}||}{||\boldsymbol{T}_0^{(s)}||}  \right)
\end{align}
where $\zeta$ is a scaling parameter. The permeability coupled model is based on the magnitude of the solid stress tensor, $||\boldsymbol{T}^{(s)}||$ normalized by the tensor magnitude of the stress in the solid due to \emph{in-situ} conditions, $||\boldsymbol{T}_0^{(s)}||$. This \emph{ad-hoc} model will serve as the basis for comparison with the coupled porosity model. In a forthcoming work, we investigate damage modified permeability models more thoroughly.

Since we are dealing with two-phase immiscible flow, we modify the 
governing equations of the \emph{quasi-static} model as follows:
\begin{subequations}
  \label{ImmiscibleMassBalanceEQ}
  \begin{align}
    \frac{\partial}{\partial t} (\rho_w^{(f)} \phi(1-S_n)) + 
    \mathrm{div}[\rho_w^{(f)}\boldsymbol{v}_w^{(f)}] &= m_w^{(f)} 
    \quad \mbox{in} \; \Omega  \\
    \frac{\partial}{\partial t} (\rho_n^{(f)} \phi S_n) + 
    \mathrm{div}[\rho_n^{(f)}\boldsymbol{v}_n^{(f)}] &= 
    m_n^{(f)} \quad \mbox{in} \; \Omega  \\
    \mathrm{div}[\boldsymbol{T}^{(s)}] &= 0 \quad 
    \mbox{in} \; \Omega \\
    p^{(f)}(\boldsymbol{x}) &=  p^\mathrm{p}(\boldsymbol{x}) 
    \quad \mbox{on} \; \Gamma^{p} \\
    \boldsymbol{u}^{(f)}(\boldsymbol{x}) &=  \boldsymbol{u}^\mathrm{p}
    (\boldsymbol{x}) \quad \mbox{on} \; \Gamma^{\boldsymbol{u}} \\
    \boldsymbol{T}^{(s)} \boldsymbol{n} &= \boldsymbol{t}^{\boldsymbol{n}} 
    \quad \mbox{on} \; \Gamma^{\boldsymbol{t}} 
  \end{align}
\end{subequations}
where $S$ is the saturation and the subscripts ``$w$'' and ``$n$'' 
denote the wetting and non-wetting phases respectively. Note that 
we have incorporated the relationship $S_n + S_w = 1.0$ in the 
first equation. The Darcy velocities for this problem are given 
as:
\begin{align}
  \boldsymbol{v}_w^{(f)} = -\alpha^{(f)}(\mathrm{grad}[p_w^{(f)}] 
  + \rho_w^{(f)} \boldsymbol{b}_w^{(f)}) \\
  \boldsymbol{v}_n^{(f)} = -\alpha^{(f)}(\mathrm{grad}[p_n^{(f)}] 
  + \rho_n^{(f)} \boldsymbol{b}_n^{(f)}) 
\end{align}
Here we have assumed no capillary pressure such that $p_n = p_w$. With this 
assumption, the equations above may be solved for $S_n$ and $p_n$ (or $p_w$). 
Note that we are using the same damage modified permeability $\alpha^{(f)}$ 
for both phases. To limit the spurious oscillations engendered by the 
advection terms, we use the CVFEM upwinding method presented by Forsyth 
in references \cite{Forsyth_1990_v5_p561,Forsyth_1991_v12_p1029}.

\begin{figure}[htb!]
  \centering
  \includegraphics[scale=1.0]{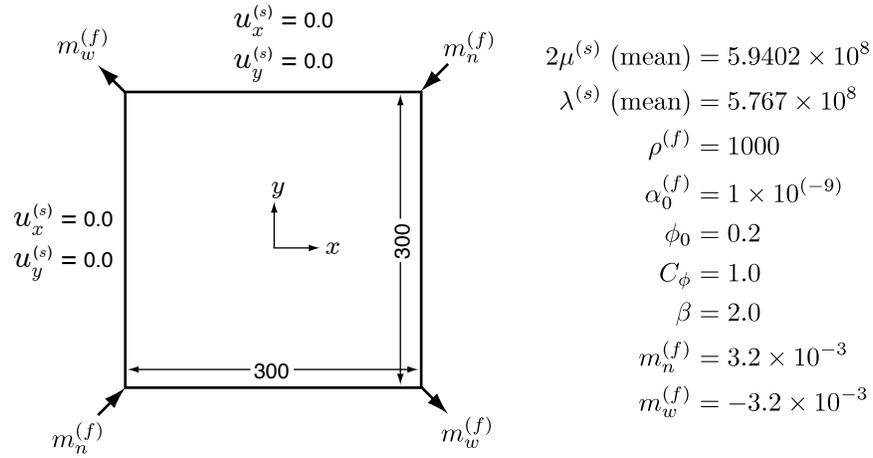}
  \caption{Two phase immiscible water flood: Problem geometry, parameters, 
  and boundary conditions.} \label{fig:FiveSpotDomain}
\end{figure}

\begin{figure}[htb!]
  \centering
\subfigure{  \includegraphics[scale=0.35]{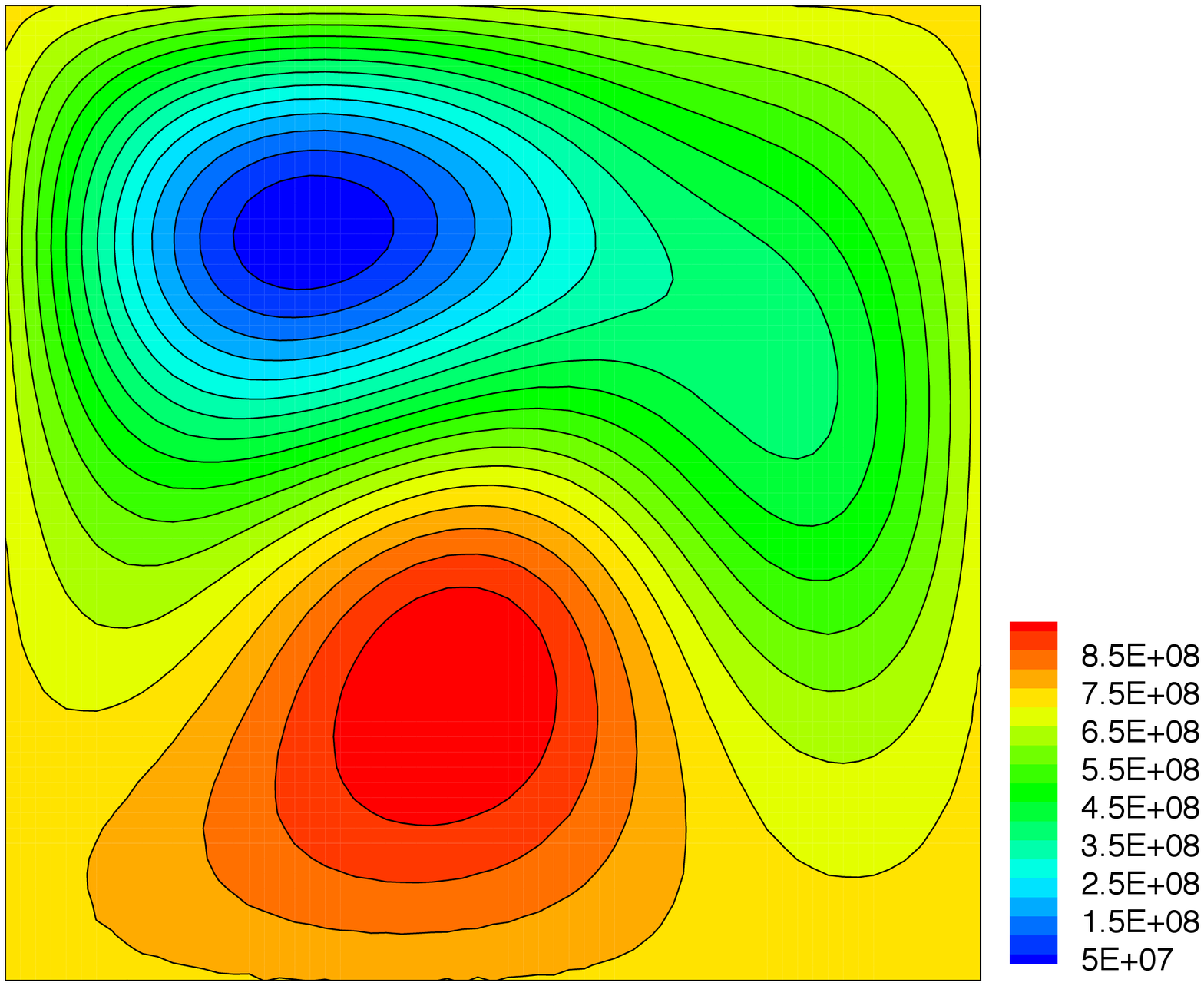}}
\subfigure{  \includegraphics[scale=0.35]{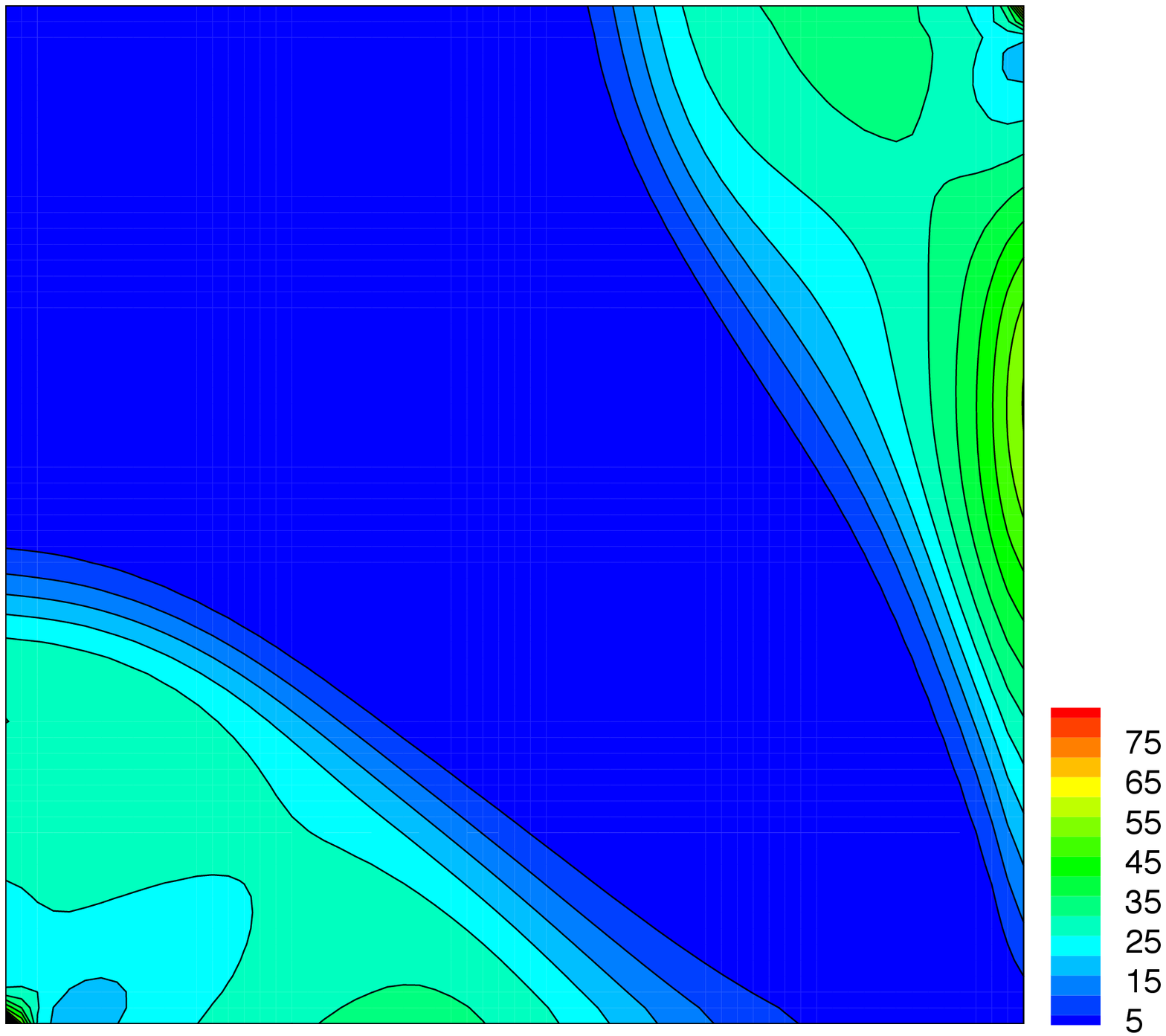}}
\subfigure{  \includegraphics[scale=0.35]{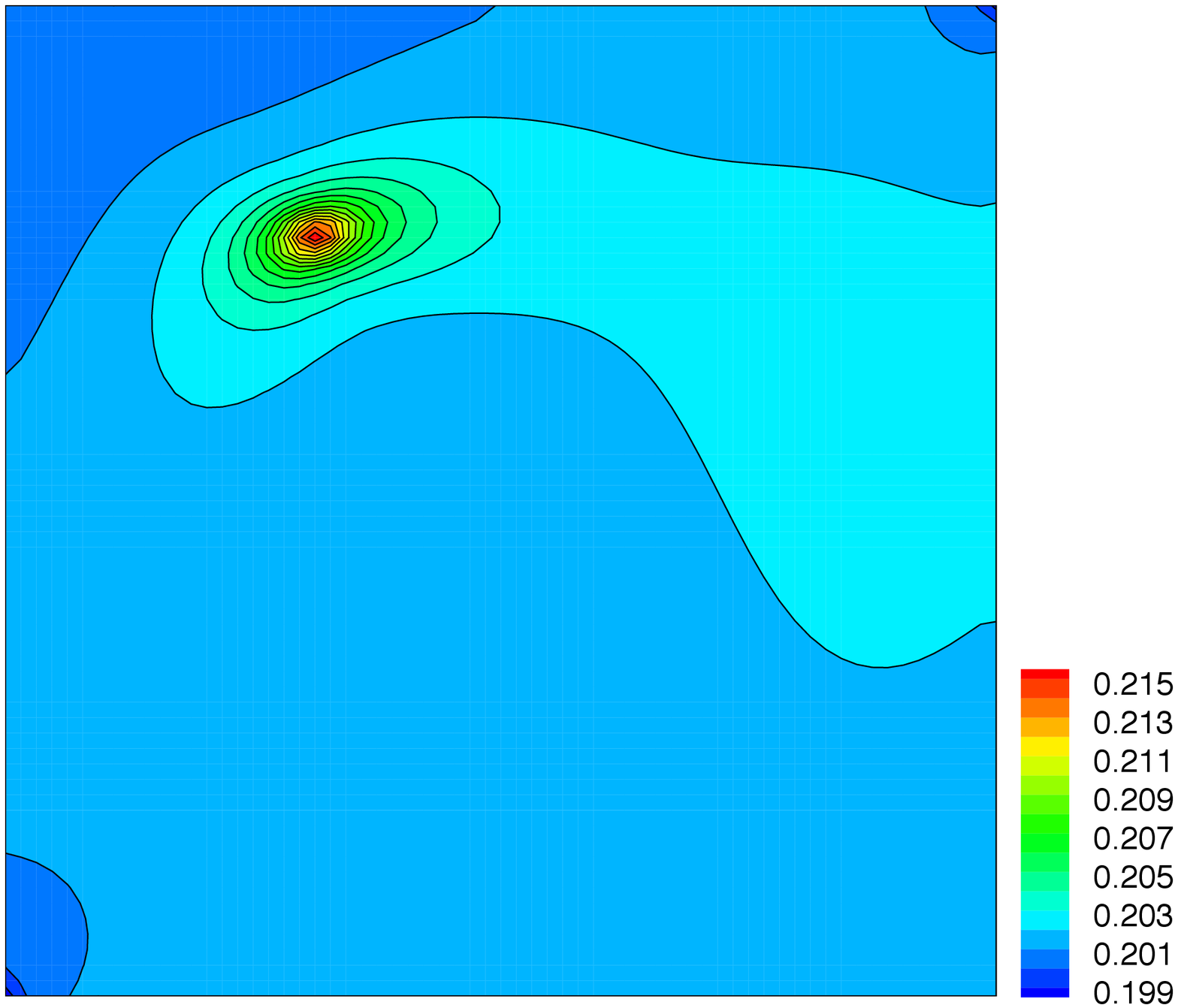}}
  \caption{Two phase immiscible water flood: Contours of $\lambda^{(s)}$ (top-left), 
  $\alpha$ (top-right), and $\phi$ (bottom).}
  \label{fig:FiveSpotLambda}
\end{figure}

\begin{figure}[htb!]
  \centering
  \includegraphics[scale=0.35]{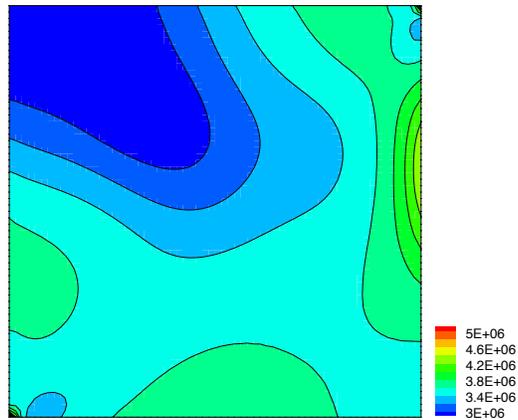}
  \caption{Two phase immiscible water flood: solid stress magnitude.}
  \label{fig:FiveSpotStress}
\end{figure}

\begin{figure}[htb!]
  \centering
\subfigure{  \includegraphics[scale=0.35]{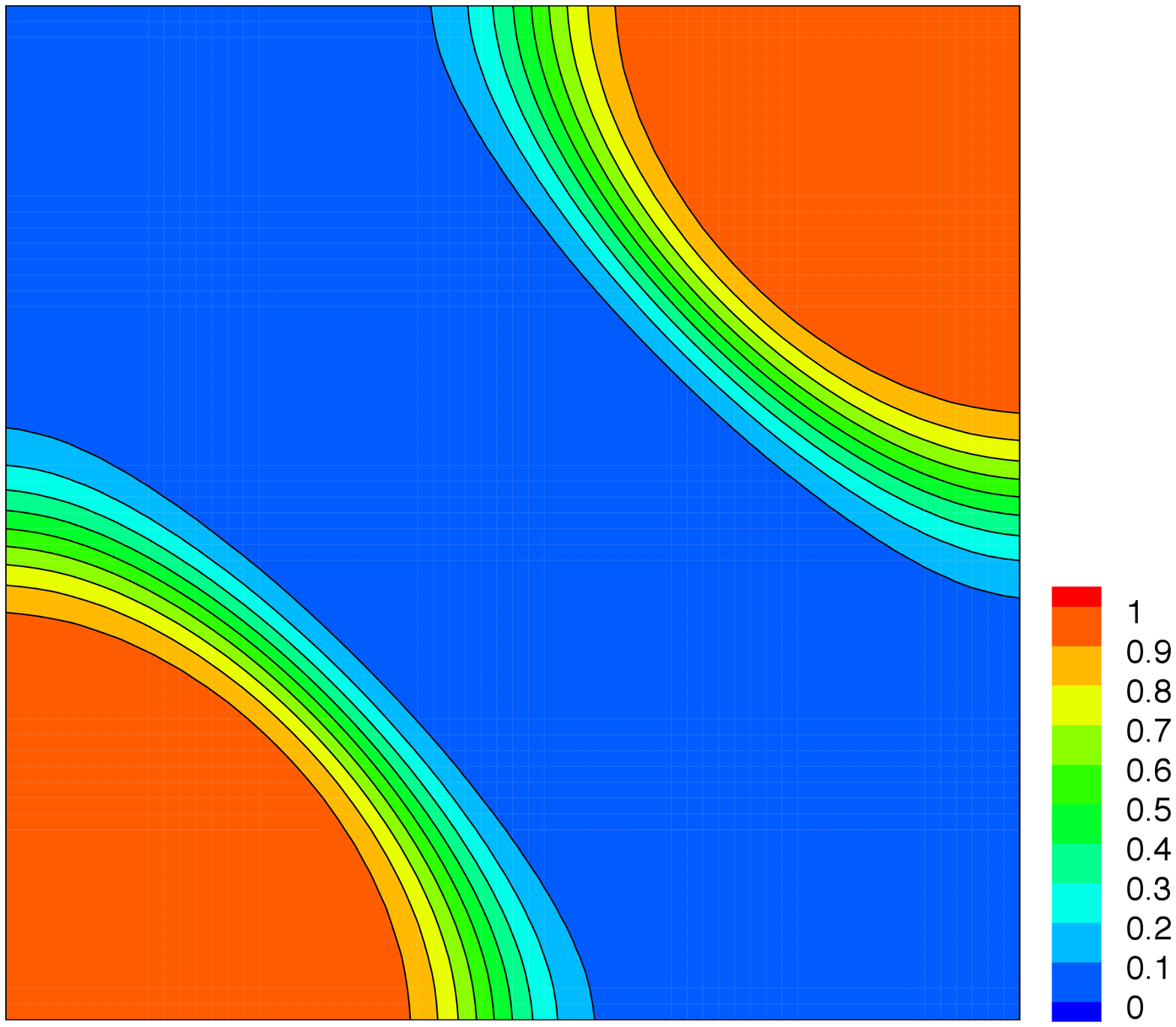}}
\subfigure{  \includegraphics[scale=0.35]{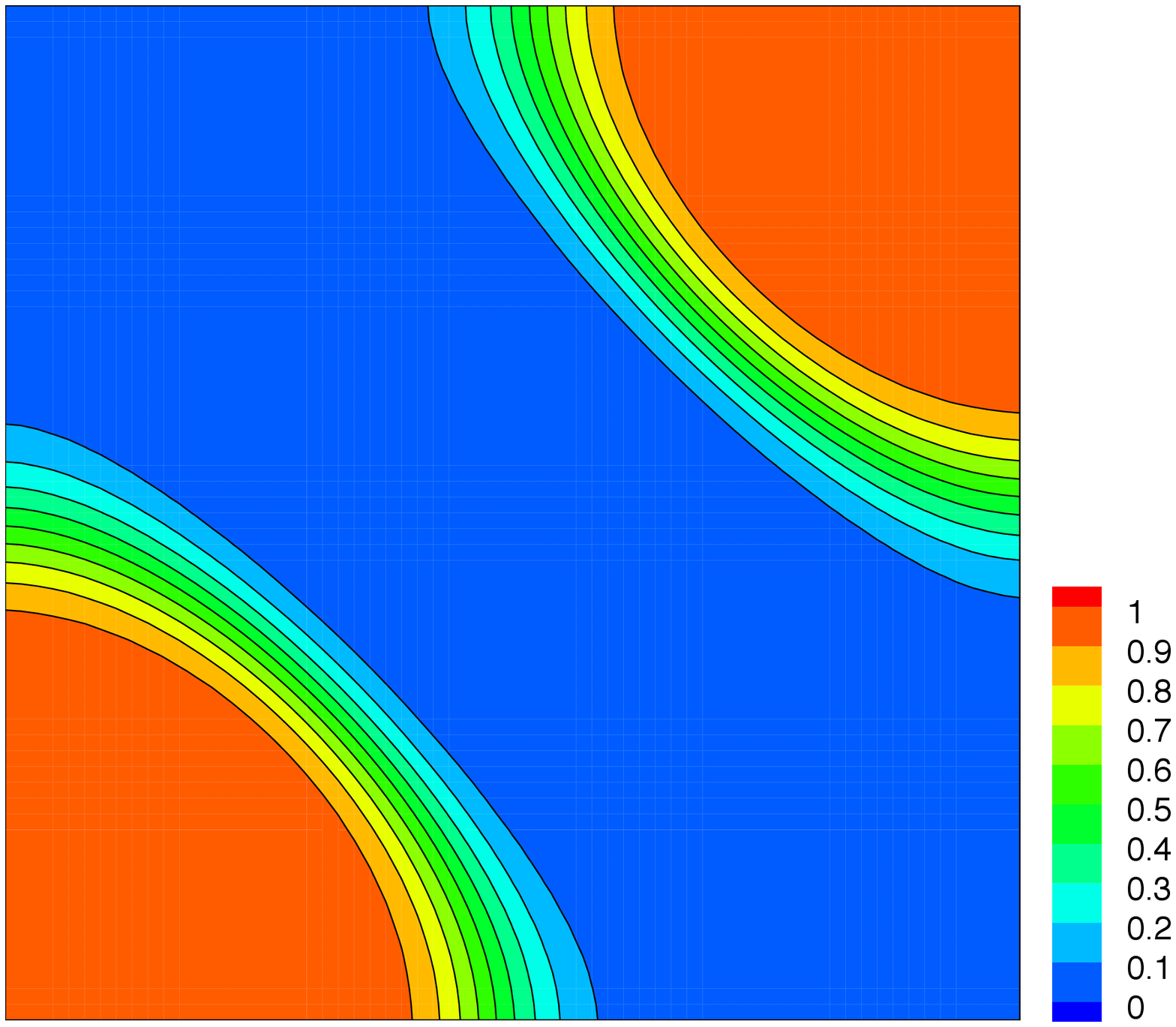}}
\subfigure{  \includegraphics[scale=0.35]{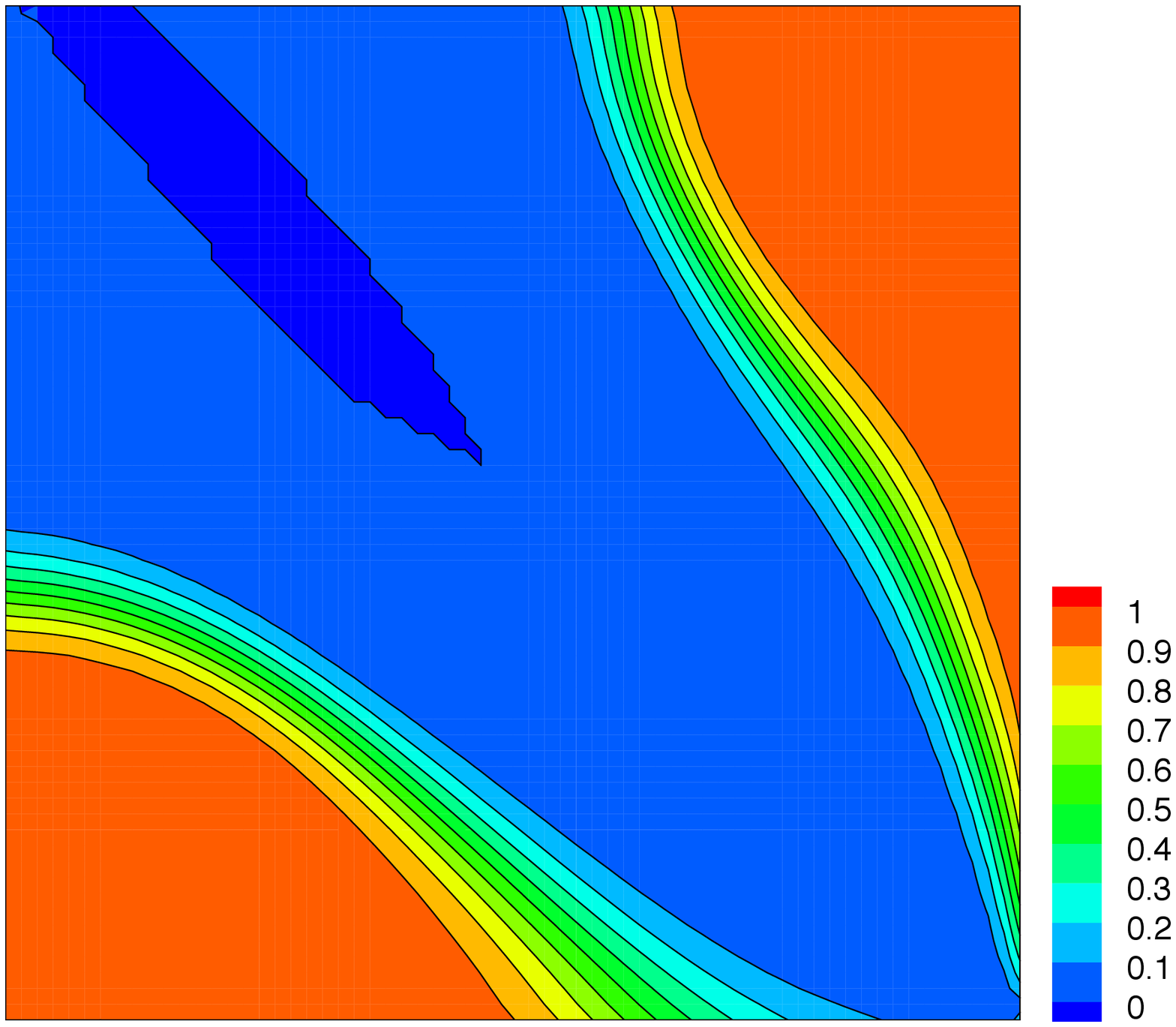}}
  \caption{Two phase immiscible water flood: Contours of $S_n$. 
  The top-left figure shows the non-wetting phase saturation for 
  uncoupled two phase model. The top-right figure shows the 
  non-wetting phase saturation for the porosity coupled model. 
  The bottom figure shows the non-wetting phase saturation 
  obtained using the damage modified permeability model. 
  There are no noticeable differences between the results 
  between the uncoupled model and the porosity coupled 
  model. In contrast, the damage modified permeability 
  model results clear show \emph{race-tracking} in the 
  regions of high stress. \label{fig:FiveSpotSaturation}}
\end{figure}

The problem domain and boundary conditions are shown in Figure \ref{fig:FiveSpotDomain}. The elastic properties for the reservoir (used only in the coupled results) were generated by the code Sierra/Encore \cite{Aria,Fuego} using a Karhunen--Lo\'eve realization. Contours of $\lambda^{(s)}$ are shown in Figure \ref{fig:FiveSpotLambda}. Notice large regions of variation of the elastic properties throughout the domain. Also shown in Figure \ref{fig:FiveSpotLambda} are contours of the porosity and permeability at time $t = 2 \times 10^{8} \; \mathrm{s}$ for both porosity coupled and permeability coupled results. Figure \ref{fig:FiveSpotStress} shows the tensor magnitude contours of the solid stress due to the pore pressure loading. When compared with Figure \ref{fig:FiveSpotLambda} it reveals that the porosity model exhibits its largest porosity in a concentrated region near the region of low solid stress magnitude. The permeability model on the other hand shows high permeability in regions of high stress as the model intends. This results in a drastically different flow behavior between the models as shown in Figure \ref{fig:FiveSpotSaturation}. 
The top-left side of Figure \ref{fig:FiveSpotSaturation} 
shows the non-wetting phase saturation for uncoupled two 
phase model. The top-right part of the figure shows the 
non-wetting phase saturation for the porosity coupled 
model. Notice that the differences with the uncoupled 
model are negligible. In contrast, the damage modified 
permeability model results clear show \emph{race-tracking} 
in the regions of high stress, which is shown the bottom 
side of Figure \ref{fig:FiveSpotSaturation}. This example 
highlights the need for further model development to 
capture damage motivated permeability effects.

\section{CONCLUDING REMARKS}
\label{Sec:Coupled_Conclusions}
We have considered a comprehensive mathematical model for 
flow of an incompressible fluid in deformable porous solid. 
The model is based on the theory of interacting continua 
with deformation-dependent porosity and damage modified 
permeability. The model also takes into account the 
dependence of viscosity the fluid on the pressure of 
the fluid. The model is fully coupled in the sense 
that the flow problem depends on the response of the 
porous solid, and the deformation of the porous solid 
depends on the velocity of the fluid. Several numerical 
coupling algorithms to obtain coupled flow-deformation 
response are also presented. Representative numerical 
results show the importance of considering the deformation 
of the solid on the flow and vice-versa. A possible future 
work is develop models and numerical coupling algorithms 
for flow in an elasto-plastic fracturable porous solid.

\section*{ACKNOWLEDGMENTS}
This work was supported in part by Sandia National Laboratories 
through the Laboratory Directed Research and Development program 
(Contract No. C11-00239). Sandia is a multiprogram laboratory 
operated by Sandia Corporation, a Lockheed Martin Company, for 
the United States Department of Energy's National Nuclear Security 
Administration under Contract DE-AC04-94AL85000.  
The second author (K.~B.~Nakshatrala) also acknowledges the 
support of the National Science Foundation under Grant no. CMMI 
1068181. The opinions expressed in this paper are those of the 
authors and do not necessarily reflect that of the sponsors.
Dr.~Martinez was supported in part by the Center for Frontiers 
of Subsurface Energy Security, an Energy Frontier Research 
Center funded by the U.S. Department of Energy, Office of 
Science, Office of Basic Energy Sciences under Award 
Number DESC0001114.

\bibliographystyle{unsrt}
\bibliography{Master_References,Books,Master_References_DAN}

\clearpage
\newpage

\end{document}